\documentclass[final]{siamltex}         

\usepackage[dvips]{color}

\usepackage{graphicx}

\usepackage{subfigure}

\graphicspath{{C:/NS/Pictures/}}

\title{Existence, uniqueness and smoothness of solution for \\3D Navier-Stokes equations with any smooth initial velocity.}

\author{ A. Tsionskiy, M. Tsionskiy  }               

\newcommand{\sz}[1]{\mbox{\Large #1 \normalsize}}
\newcommand{\szp}[1]{\mbox{\large #1 \normalsize}}

\begin{document}             

\maketitle                   

\begin{abstract}
Different authors had received a lot of results regarding the Euler and Navier-Stokes equations.
Existence and smoothness of solution for the Navier-Stokes equations in two dimensions have been known for a long time. 
Leray $\cite{jL34}$ showed that the Navier-Stokes equations in three dimensional space  have a weak solution. 
Scheffer $\cite{vS76}, \cite{vS93}$ and Shnirelman $\cite{aS97}$ obtained weak solution of the Euler equations with compact support in spacetime. Caffarelli, Kohn and Nirenberg $\cite{CKN82}$ improved Scheffer's results, and F.-H. Lin $\cite{fL98}$ simplified the proof of the results of J. Leray. Many problems and conjectures about behavior of weak solutions of the Euler and Navier-Stokes equations are described in the books of Ladyzhenskaya $\cite{oL69}$,  Bertozzi and Majda $\cite{BM02}$, Temam $\cite{RT77}$, Constantin $\cite{pC01}$ or Lemarié-Rieusset $\cite{pL02}$.

Solutions of the Navier-Stokes and Euler equations with initial conditions (Cauchy problem) for 2D and 3D cases were obtained in the converging series form by analytical iterative method using Fourier and Laplace transforms in paper $\cite{TT10}$. These solutions were received as infinitely differentiable functions. That allowed us to analyze essential aspects of the problem on a much deeper level and with more details. 
For several combinations of problem parameters numerical results were obtained and presented as graphs $\cite{TT10}$,$\;\cite{TT11}$.

This paper describes detailed proof of existence and uniqueness of the solution of the Cauchy problem for the 3D Navier- Stokes equations with any smooth initial velocity. This solution satisfies the conditions required in $\cite{CF06}$ for the problem of Navier-Stokes equations. When viscosity tends to zero this proof is correct for the Euler equations also.
\end{abstract}



\pagestyle{myheadings}
\thispagestyle{plain}
\markboth{A. TSIONSKIY, M. TSIONSKIY}{EXISTENCE, UNIQUENESS AND SMOOTHNESS OF  SOLUTION}

\section{Introduction}\ 

The solution of the Cauchy problem for the 3D Navier-Stokes equations described in this paper is based on the form of differential equations in the statement of the problem and also conditions for the initial velocity and applied force $\cite{CF06}$. It grows from classic definition of function and classic methods of analysis.


In the problem statement for Navier-Stokes equations an initial velocity is infinitely differentiable function decreasing rapidly to zero in infinity. Applied force is identically zero.

Solution of the problem will be presented by the following stages:

\textbf{First stage} (sections 2, 3) - we have moved non-linear parts of equations to the right sides.
Then we have solved the system of linear partial differential equations with constant coefficients. We have obtained the solution of this system using Fourier transforms for the space coordinates and Laplace transform for time.
From theorems about application of Fourier and Laplace transforms for system of linear partial differential equations with constant coefficients we see that in case if initial velocity and applied force are smooth enough functions decreasing in infinity, then the solution of such system is also a smooth function. (Corresponding theorems are presented in S. Bochner $\cite{SB59}$, V.P. Palamodov $\cite{VP70}$, G.E. Shilov $\cite{gS01}$, L. Hormander $\cite{LH83}$, S. Mizohata $\cite{SM73}$, J.F. Treves $\cite{JFT61}$).
Result of this stage is the integral equation for vector-function of velocity.

\textbf{Second stage} (sections 4, 5) - we have introduced perfect spaces of functions and vector-functions (I.M. Gel'fand, G.E. Chilov $\cite{GC68}$), in which we have looked 
for the solution of the problem.
We have demonstrated equivalence of the solution of the Cauchy problem in forms of differential and integral equations.

\textbf{Third stage} (section 6) - we have divided all parts of integral equation by an appropriate constant , and received
equivalent integral equation. We have also correspondingly replaced integration variables in integral operators. This newly received equivalent integral equation allowed us to analyze the Cauchy problem for the 3D Navier-Stokes equations for any value of initial fluid velocity. 

\textbf{Fourth stage} (section 6) - We have used the newly received equivalent integral equation to prove existence and uniqueness of the solution of the Cauchy problem in
all time range [0,$\infty$) based on the Caccioppoli-Banach fixed point principle (L.V. Kantorovich, G.P. Akilov $\cite{KA64},\;$ V.A. Trenogin  $\cite{VT80},\;$W. Rudin $\cite{WR73},\;$W.A. Kirk and B. Sims$\; \cite{KS01},\;$A. Granas and J. Dugundji$\; \cite{GD03},\;$J.M. Ayerbe Toledano, T. Dominguez Benavides, G. Lopez Acedo$\;\cite{ADL97}\;$).
For this purpose three following theorems were proven in this paper:

Theorem 1: Integral operator of the problem is a contraction operator;

Theorem 2: Existence and uniqueness of the solution of the problem is valid for any t $\in$ [0,$\infty$);

Theorem 3: Solution of the problem is depending on t continuously.

\textbf{Fifth stage} (section 6) - By using a priori estimation of the solution of the Cauchy problem for the 3D Navier-Stokes equations $\; \cite{oL69},\; \cite{LK63}\;$ we have shown that the energy of the whole process has a finite value for any t $\in$ [0,$\infty$).

\section{The mathematical setup}\ 

The Navier-Stokes equations describe the motion of a fluid in \(R^{N}\; ( N = 3 ) \). We look for a viscous incompressible fluid filling all of \(R^{N}\)  here. The Navier-Stokes equations are then given by

\begin{equation}\label{eqn1}
\frac{\partial u_{k}}{\partial t} \; + \; \sum_{n=1}^{N} u_{n}\frac{\partial u_{k}}{\partial x_{n}}\; =\;\nu\Delta u_{k}\; - \; \frac{\partial p}{\partial x_{k}} \; + \;f_{k}(x,t)\;\;\;\;\;(x\in R^{N},\;\;t\geq 0,\;\;{1\leq k \leq N}) 
\end{equation}
\begin{equation}\label{eqn2}
\emph{div}\,\vec{u}\;=\; \sum_{n=1}^{N} \frac{\partial u_{n}}{\partial x_{n}}\; =\;0\;\;\;\;\;\;\;\;\;\;  (x\in R^{N},t\geq 0)
\end{equation}

with initial conditions

\begin{equation}\label{eqn3}
\vec{u}(x,0)\; = \; \vec{u}^{0}(x)\;\;\;\;\;\;\;\;\;\; (x\in R^{N})
\end{equation}

Here \(\vec{u}(x,t)=(u_{k}(x,t)) \in R^{N},\;\; ({1\leq k \leq N})  \;-\; \)is an unknown velocity vector  \(( N = 3 );\; p\,(x,t)\;-\;\) is an unknown pressure; \(\vec{u}^{0}(x)\;\) is a given, \(C^{\infty}\) divergence-free vector field\(;\; f_{k}(x,t)\;\)are components of a given, externally applied force \(\vec{f}(x,t)\); \(\nu\) is a positive coefficient of the viscosity (if \(\nu = 0\) then $(\ref{eqn1})$ - $(\ref{eqn3})$ are the Euler equations); and \( \Delta\;=\; \sum_{n=1}^{N} \frac{\partial^{2}}{\partial x_{n}^{2}}\;\) is the Laplacian in the space variables. Equation $(\ref{eqn1})$ is Newton's law for a fluid element subject. Equation $(\ref{eqn2})$ says that the fluid is incompressible. For physically reasonable solutions, we accept 

\begin{equation}\label{eqn4}
u_{k}(x,t) \rightarrow 0\;\;, \;\;
\frac{\partial u_{k}}{\partial x_{n}} \;\rightarrow\;0\;\; 
\rm {as} 
\;\;\mid x \mid \;\rightarrow\; \infty\;\;\;( {1\leq k \leq N} ,\;\; {1\leq n \leq N}) \;\;\; 
\end{equation}

Hence, we will restrict attention to initial conditions $\vec{u}^{0}$ and force $\vec{f}$ that satisfy

\begin{equation}\label{eqn5}
\mid\partial_{x}^{\alpha}\vec{u}^{0}(x)\mid\;\leq\;C_{\alpha K}(1+\mid x \mid)^{-K} \quad \rm{on }\;R^{N}\;\rm{ for\; any }\;\alpha \;\rm{ and \;}K.
\end{equation}
and
\begin{equation}\label{eqn6}
\mid\partial_{x}^{\alpha}\partial_{t}^{\beta}\vec{f}(x,t)\mid\;\leq\;C_{\alpha \beta K}(1+\mid x \mid +t)^{-K} \quad \rm{on }\;R^{N}\times[0,\infty)\; \rm{ for \;any }\;\alpha,\beta \;\rm{ and \;}K.
\end{equation}

To start the process of solution let us add (\( - \sum_{n=1}^{N} u_{n}\frac{\partial u_{k}}{\partial x_{n}}\;\)) to both sides of the equations (\ref{eqn1}). Then we have:

\begin{equation}\label{eqn7}
\frac{\partial u_{k}}{\partial t}\;=\;\nu\,\Delta\,u_{k}\;-\;\frac{\partial p}{\partial x_{k}}\;+\;f_{k}(x,t)- \; \sum_{n=1}^{N} u_{n}\frac{\partial u_{k}}{\partial x_{n}}\;\;\;\;\;\;\;\;\;\;\;\; (x\in R^{N},\;\;t\geq 0,\;\;{1\leq k \leq N})
\end{equation}
\begin{equation}\label{eqn8}
\emph{div}\,\vec{u}\;=\; \sum_{n=1}^{N} \frac{\partial u_{n}}{\partial x_{n}}\; =\;0\;\;\;\;\;\;\;\;\;\;  (x\in R^{N},t\geq 0)
\end{equation}
\begin{equation}\label{eqn9}
\vec{u}(x,0)\; = \; \vec{u}^{0}(x)\;\;\;\;\;\;\;\;\;\; (x\in R^{N})
\end{equation}

\begin{equation}\label{eqn10}
u_{k}(x,t) \rightarrow 0\;\;, \;\;\frac{\partial u_{k}}{\partial x_{n}}\;\rightarrow\;0\;\; \rm{as} \;\;\mid x \mid \;\rightarrow\; \infty\;\;\;( {1\leq k \leq N} ,\;\; {1\leq n \leq N}) \;\;\; 
\end{equation}

\begin{equation}\label{eqn11}
\mid\partial_{x}^{\alpha}\vec{u}^{0}(x)\mid\;\leq\;C_{\alpha K}(1+\mid x \mid)^{-K} \quad\rm{on }\;R^{N}\;\rm{ for \; any }\;\alpha\;\rm{ and }\;K.
\end{equation}

\begin{equation}\label{eqn12}
\mid\partial_{x}^{\alpha}\partial_{t}^{\beta}\vec{f}(x,t)\mid\;\leq\;C_{\alpha \beta K}(1+\mid x \mid +t)^{-K} \quad\rm{on }\;R^{N}\times[0,\infty)\;\rm{ for \; any }\;\alpha,\beta\;\rm{ and }\;K.
\end{equation}

Let us denote

\begin{equation}\label{eqn13}
\tilde{f}_{k}(x,t)\; = \; f_{k}(x,t) \; - \; \sum_{n=1}^{N} u_{n}\frac{\partial u_{k}}{\partial x_{n}}\;\;\;\;\;\;({1\leq k \leq N})
\end{equation}

or we can present it in the vector form:

\begin{equation}\label{eqn14}
\vec{\tilde{f}}(x,t)\; = \; \vec{f}(x,t) \; - \;(\;\vec{u}\;\cdot\;\nabla\;)\;\vec{u}\;
\end{equation}
\\
\section{Solution of the system of equations $\textbf {(\ref{eqn7})}\;$ - $\;\textbf {(\ref{eqn14})}$}\ 

Let us assume that all operations below are valid. The validity of these operations will be proved in the next sections.
Taking into account our substitution (\ref{eqn13}) we see that equations $(\ref{eqn7})\; - \;(\ref{eqn9})$ are in fact system of linear partial differential equations with constant coefficients. 

Solution of this system will be presented by the following steps:

On the \textbf{first step} of our analysis we use Fourier transform $(\ref{A3})$ to solve equations $(\ref{eqn7})\; - \;(\ref{eqn14})$. 
We have got:

\[\\U_{k}( \gamma_{1}, \gamma_{2}, \gamma_{3}, t)\;=\;F[u_{k} ( x_{1}, x_{2}, x_{3}, t)]\]

\[\quad\quad\quad\quad\quad\quad -\gamma^{2}_{s}U_{k}( \gamma_{1}, \gamma_{2}, \gamma_{3}, t)\;\;=\;F\sz{[} \frac{\partial^{2}u_{k}( x_{1}, x_{2}, x_{3}, t)}{\partial x^{2}_{s}}\sz{]}\;\;\;\;\rm{[use \;(\ref{eqn10})]}\]

\[\\U_{k}^{0}( \gamma_{1} ,\gamma_{2} ,\gamma_{3})\;=\;F[u_{k}^{0} ( x_{1}, x_{2}, x_{3})]\]

\[\\P( \gamma_{1}, \gamma_{2}, \gamma_{3}, t)\;=\;F[p\, ( x_{1}, x_{2}, x_{3}, t)]\]

\[\\\tilde{F}_{k}( \gamma_{1}, \gamma_{2}, \gamma_{3}, t)\;=\;F[\tilde{f}_{k} ( x_{1}, x_{2}, x_{3}, t)]\]

\[ \\  k,s\;=\;1,2,3 \]

and then:

\begin{equation}\label{eqn134}
\frac{d U_{1}( \gamma_{1}, \gamma_{2}, \gamma_{3}, t )}{d t}  \;=\;-\nu
( \gamma_{1}^{2} +\gamma_{2}^{2} +\gamma_{3}^{2}) U_{1}( \gamma_{1}, \gamma_{2}, \gamma_{3}, t )\;+\;i\gamma_{1} P( \gamma_{1}, \gamma_{2}, \gamma_{3}, t )\;+\; \tilde{F}_{1}( \gamma_{1}, \gamma_{2}, \gamma_{3},t )
\end{equation}
\begin{equation}\label{eqn135}
\frac{d U_{2}( \gamma_{1}, \gamma_{2}, \gamma_{3}, t )}{d t}  \;=\;-\nu
( \gamma_{1}^{2} +\gamma_{2}^{2} +\gamma_{3}^{2}) U_{2}( \gamma_{1}, \gamma_{2}, \gamma_{3}, t )\;+\;i\gamma_{2} P( \gamma_{1}, \gamma_{2}, \gamma_{3}, t )\;+\; \tilde{F}_{2}( \gamma_{1}, \gamma_{2}, \gamma_{3},t )
\end{equation}
\begin{equation}\label{eqn136}
\frac{d U_{3}( \gamma_{1}, \gamma_{2}, \gamma_{3}, t )}{d t}  \;=\;-\nu
( \gamma_{1}^{2} +\gamma_{2}^{2} +\gamma_{3}^{2}) U_{3}( \gamma_{1}, \gamma_{2}, \gamma_{3}, t )\;+\;i\gamma_{3} P( \gamma_{1}, \gamma_{2}, \gamma_{3}, t )\;+\; \tilde{F}_{3}( \gamma_{1}, \gamma_{2}, \gamma_{3},t )
\end{equation}

\begin{equation}\label{eqn137}
\gamma_{1} U_{1}( \gamma_{1}, \gamma_{2}, \gamma_{3}, t ) \;+\; \gamma_{2}\, U_{2}( \gamma_{1}, \gamma_{2}, \gamma_{3}, t ) \;+\; \gamma_{3}\, U_{3}( \gamma_{1}, \gamma_{2}, \gamma_{3}, t ) \;=\;0
\end{equation}

\begin{equation}\label{eqn138}
U_{1}(\gamma_{1}, \gamma_{2},  \gamma_{3},  0)\;=\; U_{1}^{0}(\gamma_{1} ,\gamma_{2} ,\gamma_{3})
\end{equation}

\begin{equation}\label{eqn139}
U_{2}(\gamma_{1}, \gamma_{2},  \gamma_{3},  0)\;=\; U_{2}^{0}(\gamma_{1} ,\gamma_{2} ,\gamma_{3})
\end{equation}

\begin{equation}\label{eqn140}
U_{3}(\gamma_{1}, \gamma_{2},  \gamma_{3},  0)\;=\; U_{3}^{0}(\gamma_{1} ,\gamma_{2} ,\gamma_{3})
\end{equation}

Hence, we have received a system of linear ordinary differential equations with constant coefficients $(\ref{eqn134})\;-\; (\ref{eqn140})\;$ according to Fourier transforms. At the same time the initial conditions are set only for Fourier transforms of velocity components $U_{1}( \gamma_{1}, \gamma_{2}, \gamma_{3}, t ), \;U_{2}( \gamma_{1}, \gamma_{2}, \gamma_{3}, t ), \;U_{3}( \gamma_{1}, \gamma_{2}, \gamma_{3}, t )$. Because of that we can eliminate Fourier tranform for pressure $P( \gamma_{1}, \gamma_{2}, \gamma_{3}, t )$ from equations $(\ref{eqn134})\;-\; (\ref{eqn136})\;$ on the \textbf{second step} of solution.
From here assuming that $\gamma_{1} \neq 0,\;\gamma_{2} \neq 0,\;\gamma_{3} \neq 0$, we eliminate $P( \gamma_{1}, \gamma_{2}, \gamma_{3}, t )$ from equations $(\ref{eqn134})\;-\; (\ref{eqn136})\;$ and find:

\begin{eqnarray}\label{eqn141}
\frac{d}{dt} \szp{[} U_{2}( \gamma_{1}, \gamma_{2}, \gamma_{3}, t )\; -\;\frac{\gamma_{2}}{\gamma_{1}} \,U_{1}( \gamma_{1}, \gamma_{2}, \gamma_{3}, t) \szp{]} \;= -\nu( \gamma_{1}^{2} +\gamma_{2}^{2} +\gamma_{3}^{2}) \szp{[} U_{2}( \gamma_{1}, \gamma_{2}, \gamma_{3}, t )\; -\; \nonumber \\ 
\nonumber\\
-\;\frac{\gamma_{2}}{\gamma_{1}}\, U_{1}( \gamma_{1}, \gamma_{2}, \gamma_{3}, t) \szp{]}   + \; \szp{[} \tilde{F}_{2}( \gamma_{1}, \gamma_{2}, \gamma_{3}, t )\; -\;\frac{\gamma_{2}}{\gamma_{1}} \, \tilde{F}_{1}( \gamma_{1}, \gamma_{2}, \gamma_{3}, t) \szp{]} 
\quad\quad\quad\quad\quad\quad 
\end{eqnarray}

\begin{eqnarray}\label{eqn142}
\frac{d}{dt} \szp{[}  U_{3}( \gamma_{1}, \gamma_{2}, \gamma_{3}, t )\; -\;\frac{\gamma_{3}}{\gamma_{1}}\, U_{1}( \gamma_{1}, \gamma_{2}, \gamma_{3}, t) \szp{]} \;= -\nu( \gamma_{1}^{2} +\gamma_{2}^{2} +\gamma_{3}^{2}) \szp{[} U_{3}( \gamma_{1}, \gamma_{2}, \gamma_{3}, t )\; -\; \nonumber \\ 
\nonumber\\
-\;\frac{\gamma_{3}}{\gamma_{1}}\, U_{1}( \gamma_{1}, \gamma_{2}, \gamma_{3}, t) \szp{]}   + \;\szp{[} \tilde{F}_{3}( \gamma_{1}, \gamma_{2}, \gamma_{3}, t )\; -\;\frac{\gamma_{3}}{\gamma_{1}}\, \tilde{F}_{1}( \gamma_{1}, \gamma_{2}, \gamma_{3}, t) \szp{]} 
\quad\quad\quad\quad\quad\quad 
\end{eqnarray}

\begin{equation}\label{eqn143}
\gamma_{1} U_{1}( \gamma_{1}, \gamma_{2}, \gamma_{3}, t ) \;+\; \gamma_{2}\, U_{2}( \gamma_{1}, \gamma_{2}, \gamma_{3}, t ) \;+\; \gamma_{3}\, U_{3}( \gamma_{1}, \gamma_{2}, \gamma_{3}, t ) \;=\;0
\end{equation}

\begin{equation}\label{eqn144}
U_{1}(\gamma_{1}, \gamma_{2},  \gamma_{3},  0)\;=\; U_{1}^{0}(\gamma_{1} ,\gamma_{2} ,\gamma_{3})
\end{equation}

\begin{equation}\label{eqn145}
U_{2}(\gamma_{1}, \gamma_{2},  \gamma_{3},  0)\;=\; U_{2}^{0}(\gamma_{1} ,\gamma_{2} ,\gamma_{3})
\end{equation}

\begin{equation}\label{eqn146}
U_{3}(\gamma_{1}, \gamma_{2},  \gamma_{3},  0)\;=\; U_{3}^{0}(\gamma_{1} ,\gamma_{2} ,\gamma_{3})
\end{equation}

On the \textbf{third step} we use Laplace transform $(\ref{A4}), (\ref{A5})$  for a system of linear ordinary differential equations with constant coefficients $(\ref{eqn141})\;-\; (\ref{eqn143})\;$  and have as a result the system of linear algebraic equations with constant coefficients:

\[U_{k}^{\otimes} (\gamma_{1}, \gamma_{2}, \gamma_{3}, \eta) \;=\;L[U_{k}(\gamma_{1}, \gamma_{2}, \gamma_{3}, t)] \;\;\;\;\;\;\;      \rm{k=1,2,3}\] 
\[\tilde{F}_{k}^{\otimes} (\gamma_{1}, \gamma_{2}, \gamma_{3}, \eta) \;=\;L[\tilde{F}_{k}(\gamma_{1}, \gamma_{2}, \gamma_{3}, t)] \;\;\;\;\;\;\;      \rm{k=1,2,3}\] 

\begin{eqnarray}\label{eqn147}
\eta \szp{[}  U_{2}^{\otimes}( \gamma_{1}, \gamma_{2}, \gamma_{3}, \eta )\; -\;\frac{\gamma_{2}}{\gamma_{1}} U_{1}^{\otimes}( \gamma_{1}, \gamma_{2}, \gamma_{3}, \eta) \szp{]}  \;-\; \szp{[}   U_{2}( \gamma_{1}, \gamma_{2}, \gamma_{3}, 0 )\; -\;\frac{\gamma_{2}}{\gamma_{1}} U_{1}( \gamma_{1}, \gamma_{2}, \gamma_{3}, 0) \szp{]}  \;=  \nonumber
\\
\nonumber\\
-\nu( \gamma_{1}^{2} +\gamma_{2}^{2} +\gamma_{3}^{2})\szp{[}   U_{2}^{\otimes}                          ( \gamma_{1}, \gamma_{2}, \gamma_{3}, \eta )\; -\;\frac{\gamma_{2}}{\gamma_{1}} U_{1}^{\otimes}( \gamma_{1}, \gamma_{2}, \gamma_{3}, \eta) \szp{]}  \;+ 
\quad\quad\quad\quad\quad\quad 
\nonumber
\\
\nonumber\\
+\; \szp{[}   \tilde{F}_{2}^{\otimes}( \gamma_{1}, \gamma_{2}, \gamma_{3}, \eta )\; -\;\frac{\gamma_{2}}{\gamma_{1}} \tilde{F}_{1}^{\otimes}( \gamma_{1}, \gamma_{2}, \gamma_{3}, \eta) \szp{]}
\quad\quad\quad\quad\quad\quad \quad\quad\quad \quad\quad\quad 
\end{eqnarray}

\begin{eqnarray}\label{eqn148}
\eta  \szp{[}   U_{3}^{\otimes}( \gamma_{1}, \gamma_{2}, \gamma_{3}, \eta )\; -\;\frac{\gamma_{3}}{\gamma_{1}} U_{1}^{\otimes}( \gamma_{1}, \gamma_{2}, \gamma_{3}, \eta) \szp{]}   \;-\; \szp{[}    U_{3}( \gamma_{1}, \gamma_{2}, \gamma_{3}, 0 )\; -\;\frac{\gamma_{3}}{\gamma_{1}} U_{1}( \gamma_{1}, \gamma_{2}, \gamma_{3}, 0) \szp{]}   \;=  \nonumber
\\
\nonumber\\
-\nu( \gamma_{1}^{2} +\gamma_{2}^{2} +\gamma_{3}^{2})\szp{[}    U_{3}^{\otimes}                          ( \gamma_{1}, \gamma_{2}, \gamma_{3}, \eta )\; -\;\frac{\gamma_{3}}{\gamma_{1}} U_{1}^{\otimes}( \gamma_{1}, \gamma_{2}, \gamma_{3}, \eta) \szp{]}   \;+ 
\quad\quad\quad\quad\quad\quad 
\nonumber
\\
\nonumber\\
+\; \szp{[}    \tilde{F}_{3}^{\otimes}( \gamma_{1}, \gamma_{2}, \gamma_{3}, \eta )\; -\;\frac{\gamma_{3}}{\gamma_{1}} \tilde{F}_{1}^{\otimes}( \gamma_{1}, \gamma_{2}, \gamma_{3}, \eta) \szp{]} 
\quad\quad\quad\quad\quad\quad \quad\quad\quad \quad\quad\quad 
\end{eqnarray}

\begin{equation}\label{eqn149}
\gamma_{1} U_{1}^{\otimes}( \gamma_{1}, \gamma_{2}, \gamma_{3}, \eta ) \;+\; \gamma_{2}\, U_{2}^{\otimes}( \gamma_{1}, \gamma_{2}, \gamma_{3}, \eta ) \;+\; \gamma_{3}\, U_{3}^{\otimes}( \gamma_{1}, \gamma_{2}, \gamma_{3}, \eta ) \;=\;0
\end{equation}

\begin{equation}\label{eqn150}
U_{1}(\gamma_{1}, \gamma_{2},  \gamma_{3},  0)\;=\; U_{1}^{0}(\gamma_{1} ,\gamma_{2} ,\gamma_{3})
\end{equation}

\begin{equation}\label{eqn151}
U_{2}(\gamma_{1}, \gamma_{2},  \gamma_{3},  0)\;=\; U_{2}^{0}(\gamma_{1} ,\gamma_{2} ,\gamma_{3})
\end{equation}

\begin{equation}\label{eqn152}
U_{3}(\gamma_{1}, \gamma_{2},  \gamma_{3},  0)\;=\; U_{3}^{0}(\gamma_{1} ,\gamma_{2} ,\gamma_{3})
\end{equation}

Let us rewrite system of equations $(\ref{eqn147})\;-\; (\ref{eqn149})\;$ in the following form:

\begin{eqnarray}\label{eqn147a}
\szp{[}\eta \;+\;\nu( \gamma_{1}^{2} +\gamma_{2}^{2} +\gamma_{3}^{2})\szp{]}\frac{\gamma_{2}}{\gamma_{1}} U_{1}^{\otimes}( \gamma_{1}, \gamma_{2}, \gamma_{3}, \eta)   \;-\; \szp{[}\eta \;+\;\nu( \gamma_{1}^{2} +\gamma_{2}^{2} +\gamma_{3}^{2})\szp{]} U_{2}^{\otimes}( \gamma_{1}, \gamma_{2}, \gamma_{3}, \eta)  \;=  \nonumber
\\
\nonumber\\
\;\;\; \szp{[}  \frac{\gamma_{2}}{\gamma_{1}} \tilde{F}_{1}^{\otimes}( \gamma_{1}, \gamma_{2}, \gamma_{3}, \eta) \; -\;\tilde{F}_{2}^{\otimes}( \gamma_{1}, \gamma_{2}, \gamma_{3}, \eta ) \szp{]}+ \szp{[} \frac{\gamma_{2}}{\gamma_{1}} U_{1}( \gamma_{1}, \gamma_{2}, \gamma_{3}, 0)  \; -\;U_{2}( \gamma_{1}, \gamma_{2}, \gamma_{3}, 0 ) \szp{]}
\quad
\end{eqnarray}

\begin{eqnarray}\label{eqn148a}
\szp{[}\eta \;+\;\nu( \gamma_{1}^{2} +\gamma_{2}^{2} +\gamma_{3}^{2})\szp{]}\frac{\gamma_{3}}{\gamma_{1}} U_{1}^{\otimes}( \gamma_{1}, \gamma_{2}, \gamma_{3}, \eta)   \;-\; \szp{[}\eta \;+\;\nu( \gamma_{1}^{2} +\gamma_{2}^{2} +\gamma_{3}^{2})\szp{]} U_{3}^{\otimes}( \gamma_{1}, \gamma_{2}, \gamma_{3}, \eta)  \;=  \nonumber
\\
\nonumber\\
\;\;\; \szp{[}  \frac{\gamma_{3}}{\gamma_{1}} \tilde{F}_{1}^{\otimes}( \gamma_{1}, \gamma_{2}, \gamma_{3}, \eta) \; -\;\tilde{F}_{3}^{\otimes}( \gamma_{1}, \gamma_{2}, \gamma_{3}, \eta ) \szp{]}+ \szp{[} \frac{\gamma_{3}}{\gamma_{1}} U_{1}( \gamma_{1}, \gamma_{2}, \gamma_{3}, 0)  \; -\;U_{3}( \gamma_{1}, \gamma_{2}, \gamma_{3}, 0 ) \szp{]}
\quad
\end{eqnarray}

\begin{equation}\label{eqn149a}
\gamma_{1} U_{1}^{\otimes}( \gamma_{1}, \gamma_{2}, \gamma_{3}, \eta ) \;+\; \gamma_{2}\, U_{2}^{\otimes}( \gamma_{1}, \gamma_{2}, \gamma_{3}, \eta ) \;+\; \gamma_{3}\, U_{3}^{\otimes}( \gamma_{1}, \gamma_{2}, \gamma_{3}, \eta ) \;=\;0
\end{equation}

Determinant of this system is

\begin{eqnarray}\label{eqn1481b}
\Delta = \left| \begin{array}{ccc}
\szp{[}\eta \;+\;\nu( \gamma_{1}^{2} +\gamma_{2}^{2} +\gamma_{3}^{2})\szp{]}\frac{\gamma_{2}}{\gamma_{1}} & -\szp{[}\eta \;+\;\nu( \gamma_{1}^{2} +\gamma_{2}^{2} +\gamma_{3}^{2})\szp{]} & 0 \\
\nonumber\\
 \szp{[}\eta \;+\;\nu( \gamma_{1}^{2} +\gamma_{2}^{2} +\gamma_{3}^{2})\szp{]}\frac{\gamma_{3}}{\gamma_{1}} & 0 & -\szp{[}\eta \;+\;\nu( \gamma_{1}^{2} +\gamma_{2}^{2} +\gamma_{3}^{2})\szp{]} \\
\nonumber\\ 
\gamma_{1} & \gamma_{2} & \gamma_{3} \end{array} \right| = 
\nonumber\\ 
\nonumber\\
\nonumber\\
=\;\frac{\szp{[}\eta \;+\;\nu( \gamma_{1}^{2} +\gamma_{2}^{2} +\gamma_{3}^{2})\szp{]}^{2}( \gamma_{1}^{2} +\gamma_{2}^{2} +\gamma_{3}^{2})}{\gamma_{1}}
\neq \;\; 0
\quad\quad\quad\quad\quad\quad\quad\quad\quad\quad\quad\quad
\end{eqnarray}

And consequently the system of equations $(\ref{eqn147})\;-\; (\ref{eqn149})\;$ and/or $\;(\ref{eqn147a})\;-\; (\ref{eqn149a})\;$ has a unique solution. Taking into account formulas $(\ref{eqn150})\;-\; (\ref{eqn152})\;$ we can write this solution in the following form:

\begin{eqnarray}\label{eqn153}
U_{1}^{\otimes}( \gamma_{1}, \gamma_{2}, \gamma_{3}, \eta )\;=\;\frac{[( \gamma_{2}^{2} +\gamma_{3}^{2})  \tilde{F}_{1}^{\otimes}( \gamma_{1}, \gamma_{2}, \gamma_{3}, \eta) - \gamma_{1}\gamma_{2} \tilde{F}_{2}^{\otimes}( \gamma_{1}, \gamma_{2}, \gamma_{3}, \eta) - \gamma_{1}\gamma_{3} \tilde{F}_{3}^{\otimes}( \gamma_{1}, \gamma_{2}, \gamma_{3}, \eta)]}{ (\gamma_{1}^{2} +\gamma_{2}^{2} +\gamma_{3}^{2}) [\eta+\nu (\gamma_{1}^{2} +\gamma_{2}^{2} +\gamma_{3}^{2})] }\;+\nonumber
\\
\nonumber\\
+\; \frac{ U_{1}^{0}(\gamma_{1} , \gamma_{2} , \gamma_{3})}{[\eta+\nu (\gamma_{1}^{2} +\gamma_{2}^{2} +\gamma_{3}^{2})] } 
\quad\quad\quad\quad\quad\quad \quad\quad\quad \quad\quad\quad 
\quad\quad\quad\quad\quad\quad 
\end{eqnarray}

\begin{eqnarray}\label{eqn154}
U_{2}^{\otimes}( \gamma_{1}, \gamma_{2}, \gamma_{3}, \eta )\;=\;\frac{[( \gamma_{3}^{2} +\gamma_{1}^{2})  \tilde{F}_{2}^{\otimes}( \gamma_{1}, \gamma_{2}, \gamma_{3}, \eta) - \gamma_{2}\gamma_{3} \tilde{F}_{3}^{\otimes}( \gamma_{1}, \gamma_{2}, \gamma_{3}, \eta) - \gamma_{2}\gamma_{1} \tilde{F}_{1}^{\otimes}( \gamma_{1}, \gamma_{2}, \gamma_{3}, \eta)]}{ (\gamma_{1}^{2} +\gamma_{2}^{2} +\gamma_{3}^{2}) [\eta+\nu (\gamma_{1}^{2} +\gamma_{2}^{2} +\gamma_{3}^{2})] }\;+\nonumber
\\
\nonumber\\
+\; \frac{ U_{2}^{0}(\gamma_{1} , \gamma_{2} , \gamma_{3})}{[\eta+\nu (\gamma_{1}^{2} +\gamma_{2}^{2} +\gamma_{3}^{2})] } 
\quad\quad\quad\quad\quad\quad \quad\quad\quad \quad\quad\quad 
\quad\quad\quad\quad\quad\quad 
\end{eqnarray}

\begin{eqnarray}\label{eqn155}
U_{3}^{\otimes}( \gamma_{1}, \gamma_{2}, \gamma_{3}, \eta )\;=\;\frac{[( \gamma_{1}^{2} +\gamma_{2}^{2})  \tilde{F}_{3}^{\otimes}( \gamma_{1}, \gamma_{2}, \gamma_{3}, \eta) - \gamma_{3}\gamma_{1} \tilde{F}_{1}^{\otimes}( \gamma_{1}, \gamma_{2}, \gamma_{3}, \eta) - \gamma_{3}\gamma_{2} \tilde{F}_{2}^{\otimes}( \gamma_{1}, \gamma_{2}, \gamma_{3}, \eta)]}{ (\gamma_{1}^{2} +\gamma_{2}^{2} +\gamma_{3}^{2}) [\eta+\nu (\gamma_{1}^{2} +\gamma_{2}^{2} +\gamma_{3}^{2})] }\;+ \nonumber
\\
\nonumber\\
+\; \frac{ U_{3}^{0}(\gamma_{1} , \gamma_{2} , \gamma_{3})}{[\eta+\nu (\gamma_{1}^{2} +\gamma_{2}^{2} +\gamma_{3}^{2})] } 
\quad\quad\quad\quad\quad\quad \quad\quad\quad \quad\quad\quad 
\quad\quad\quad\quad\quad\quad 
\end{eqnarray}

Then we use the convolution theorem with the convolution formula (\ref{A6}) and integral (\ref{A7}) for $(\ref{eqn153})\;-\; (\ref{eqn155})\;$and obtain:

\begin{eqnarray}\label{eqn156}
U_{1}(\gamma_{1}, \gamma_{2}, \gamma_{3}, t)\;=\;
\quad\quad\quad\quad\quad\quad \quad\quad\quad \quad\quad\quad 
\quad\quad\quad\quad\quad\quad
 \nonumber
\\
\nonumber\\
\int_{0}^{t} \sz{e} ^{-\nu (\gamma_{1}^{2} +\gamma_{2}^{2} +\gamma_{3}^{2}) (t-\tau)} \frac{[( \gamma_{2}^{2} +\gamma_{3}^{2})  \tilde{F}_{1}( \gamma_{1}, \gamma_{2}, \gamma_{3}, \tau) -\gamma_{1}\gamma_{2} \tilde{F}_{2}( \gamma_{1}, \gamma_{2}, \gamma_{3}, \tau ) -\gamma_{1}\gamma_{3} \tilde{F}_{3}( \gamma_{1}, \gamma_{2}, \gamma_{3}, \tau ) ]}{ (\gamma_{1}^{2} +\gamma_{2}^{2}+\gamma_{3}^{2} ) }\,d\tau\;+ \nonumber
\\
\nonumber\\
+\;\sz{e} ^{-\nu (\gamma_{1}^{2} +\gamma_{2}^{2} +\gamma_{3}^{2}) t} \;U_{1}^{0}(\gamma_{1} ,\gamma_{2} ,\gamma_{3})
\quad\quad\quad\quad\quad\quad \quad\quad\quad \quad\quad\quad 
\quad\quad\quad\quad 
\end{eqnarray}

\begin{eqnarray}\label{eqn157}
U_{2}(\gamma_{1}, \gamma_{2}, \gamma_{3}, t)\;=\; 
\quad\quad\quad\quad\quad\quad \quad\quad\quad \quad\quad\quad 
\quad\quad\quad\quad\quad\quad
\nonumber
\\
\nonumber\\
\int_{0}^{t} \sz{e} ^{-\nu (\gamma_{1}^{2} +\gamma_{2}^{2} +\gamma_{3}^{2}) (t-\tau)} \frac{[( \gamma_{3}^{2} +\gamma_{1}^{2})  \tilde{F}_{2}( \gamma_{1}, \gamma_{2}, \gamma_{3}, \tau) -\gamma_{2}\gamma_{3} \tilde{F}_{3}( \gamma_{1}, \gamma_{2}, \gamma_{3}, \tau ) -\gamma_{2}\gamma_{1} \tilde{F}_{1}( \gamma_{1}, \gamma_{2}, \gamma_{3}, \tau ) ]}{ (\gamma_{1}^{2} +\gamma_{2}^{2}+\gamma_{3}^{2} ) }\,d\tau\;+ \nonumber
\\
\nonumber\\
+\;\sz{e} ^{-\nu (\gamma_{1}^{2} +\gamma_{2}^{2} +\gamma_{3}^{2}) t} \;U_{2}^{0}(\gamma_{1} ,\gamma_{2} ,\gamma_{3})
\quad\quad\quad\quad\quad\quad \quad\quad\quad \quad\quad\quad 
\quad\quad\quad\quad 
\end{eqnarray}

\begin{eqnarray}\label{eqn158}
U_{3}(\gamma_{1}, \gamma_{2}, \gamma_{3}, t)\;=\; 
\quad\quad\quad\quad\quad\quad \quad\quad\quad \quad\quad\quad 
\quad\quad\quad\quad\quad\quad
\nonumber
\\
\nonumber\\
\int_{0}^{t} \sz{e} ^{-\nu (\gamma_{1}^{2} +\gamma_{2}^{2} +\gamma_{3}^{2}) (t-\tau)} \frac{[( \gamma_{1}^{2} +\gamma_{2}^{2})  \tilde{F}_{3}( \gamma_{1}, \gamma_{2}, \gamma_{3}, \tau) -\gamma_{3}\gamma_{1} \tilde{F}_{1}( \gamma_{1}, \gamma_{2}, \gamma_{3}, \tau ) -\gamma_{3}\gamma_{2} \tilde{F}_{2}( \gamma_{1}, \gamma_{2}, \gamma_{3}, \tau ) ]}{ (\gamma_{1}^{2} +\gamma_{2}^{2}+\gamma_{3}^{2} ) }\,d\tau\;+ 
\nonumber
\\
\nonumber\\
+\;\sz{e} ^{-\nu (\gamma_{1}^{2} +\gamma_{2}^{2} +\gamma_{3}^{2}) t} \;U_{3}^{0}(\gamma_{1} ,\gamma_{2} ,\gamma_{3})
\quad\quad\quad\quad\quad\quad \quad\quad\quad \quad\quad\quad 
\quad\quad\quad\quad 
\end{eqnarray}
\\
Using the Fourier inversion formula $(\ref{A3})$ we receive:

\begin{eqnarray}\label{eqn160}
u_{1}(x_{1}, x_{2}, x_{3}, t)\;=\; 
\frac{1}{(2\pi)^{3/2}} \int_{-\infty}^{\infty} \int_{-\infty}^{\infty} \int_{-\infty}^{\infty} \biggl[ \int_{0}^{t} \sz{e} ^{-\nu (\gamma_{1}^{2} +\gamma_{2}^{2} +\gamma_{3}^{2}) (t-\tau)} \frac{[ ( \gamma_{2}^{2} +\gamma_{3}^{2})  \tilde{F}_{1}( \gamma_{1}, \gamma_{2}, \gamma_{3}, \tau)]} { (\gamma_{1}^{2} +\gamma_{2}^{2}+\gamma_{3}^{2} ) } \,d\tau \;- \nonumber
\\
\nonumber\\
-\; \int_{0}^{t} \sz{e} ^{-\nu (\gamma_{1}^{2} +\gamma_{2}^{2} +\gamma_{3}^{2}) (t-\tau)}\frac{ [\gamma_{1}\gamma_{2} \tilde{F}_{2}( \gamma_{1}, \gamma_{2}, \gamma_{3}, \tau ) + \gamma_{1}\gamma_{3} \tilde{F}_{3}( \gamma_{1}, \gamma_{2}, \gamma_{3}, \tau )]} { (\gamma_{1}^{2} +\gamma_{2}^{2}+\gamma_{3}^{2} ) }\,d\tau\;   +
\quad
\nonumber\\
\nonumber\\
+\;\sz{e} ^{-\nu (\gamma_{1}^{2} +\gamma_{2}^{2} +\gamma_{3}^{2}) t} \;U_{1}^{0}(\gamma_{1} ,\gamma_{2} ,\gamma_{3})\biggr] \;\sz{e} ^{-i(x_{1}\gamma_{1}+x_{2}\gamma_{2}+x_{3}\gamma_{3})}\,d\gamma_{1}d\gamma_{2}d\gamma_{3}\;=
\nonumber\\
\nonumber\\
=\;\frac{1}{8\pi^{3}} \int_{-\infty}^{\infty} \int_{-\infty}^{\infty} \int_{-\infty}^{\infty}\frac{( \gamma_{2}^{2} +\gamma_{3}^{2})} { (\gamma_{1}^{2} +\gamma_{2}^{2}+\gamma_{3}^{2} ) }  \biggl[ \int_{0}^{t} \sz{e} ^{-\nu (\gamma_{1}^{2} +\gamma_{2}^{2} +\gamma_{3}^{2}) (t-\tau)} \int_{-\infty}^{\infty}\int_{-\infty}^{\infty}\int_{-\infty}^{\infty}\sz{e}  ^{i(\tilde x_{1}\gamma_{1}+\tilde x_{2}\gamma_{2}+\tilde x_{3}\gamma_{3})} \cdot 
\nonumber\\
\nonumber\\
\cdot \tilde{f}_{1}(\tilde x_{1},\tilde x_{2}, \tilde x_{3},\tau)\,d\tilde x_{1}d\tilde x_{2} d\tilde x_{3}d\tau\biggr]\sz{e} ^{-i(x_{1}\gamma_{1}+x_{2}\gamma_{2}+x_{3}\gamma_{3})}\,d\gamma_{1}d\gamma_{2}d\gamma_{3}\;-
\nonumber\\
\nonumber\\
-\;\frac{1}{8\pi^{3}} \int_{-\infty}^{\infty} \int_{-\infty}^{\infty} \int_{-\infty}^{\infty}\frac{ \gamma_{1}\gamma_{2}} { (\gamma_{1}^{2} +\gamma_{2}^{2}+\gamma_{3}^{2} ) }  \biggl[ \int_{0}^{t} \sz{e} ^{-\nu (\gamma_{1}^{2} +\gamma_{2}^{2} +\gamma_{3}^{2}) (t-\tau)} \int_{-\infty}^{\infty}\int_{-\infty}^{\infty}\int_{-\infty}^{\infty}\sz{e}  ^{i(\tilde x_{1}\gamma_{1}+\tilde x_{2}\gamma_{2}+\tilde x_{3}\gamma_{3})} \cdot 
\nonumber\\
\nonumber\\
\cdot \tilde{f}_{2}(\tilde x_{1},\tilde x_{2}, \tilde x_{3},\tau)\,d\tilde x_{1}d\tilde x_{2} d\tilde x_{3}d\tau\biggr]\sz{e} ^{-i(x_{1}\gamma_{1}+x_{2}\gamma_{2}+x_{3}\gamma_{3})}\,d\gamma_{1}d\gamma_{2}d\gamma_{3}\;-
\nonumber\\
\nonumber\\
-\;\frac{1}{8\pi^{3}} \int_{-\infty}^{\infty} \int_{-\infty}^{\infty} \int_{-\infty}^{\infty}\frac{ \gamma_{1}\gamma_{3}} { (\gamma_{1}^{2} +\gamma_{2}^{2}+\gamma_{3}^{2} ) }  \biggl[ \int_{0}^{t} \sz{e} ^{-\nu (\gamma_{1}^{2} +\gamma_{2}^{2} +\gamma_{3}^{2}) (t-\tau)} \int_{-\infty}^{\infty}\int_{-\infty}^{\infty}\int_{-\infty}^{\infty}\sz{e}  ^{i(\tilde x_{1}\gamma_{1}+\tilde x_{2}\gamma_{2}+\tilde x_{3}\gamma_{3})} \cdot 
\nonumber\\
\nonumber\\
\cdot \tilde{f}_{3}(\tilde x_{1},\tilde x_{2}, \tilde x_{3},\tau)\,d\tilde x_{1}d\tilde x_{2} d\tilde x_{3}d\tau\biggr]\sz{e} ^{-i(x_{1}\gamma_{1}+x_{2}\gamma_{2}+x_{3}\gamma_{3})}\,d\gamma_{1}d\gamma_{2}d\gamma_{3}\;+
\nonumber\\
\nonumber\\
+\;\frac{1}{8\pi^{3}} \int_{-\infty}^{\infty} \int_{-\infty}^{\infty} \int_{-\infty}^{\infty} \sz{e} ^{-\nu (\gamma_{1}^{2} +\gamma_{2}^{2} +\gamma_{3}^{2}) t}\biggl[  \int_{-\infty}^{\infty}\int_{-\infty}^{\infty}\int_{-\infty}^{\infty}\sz{e}  ^{i(\tilde x_{1}\gamma_{1}+\tilde x_{2}\gamma_{2}+\tilde x_{3}\gamma_{3})}
 \cdot \quad\quad\quad\quad\quad\quad\quad\quad\quad
\nonumber\\
\nonumber\\
\cdot \; u_{1}^{0}(\tilde x_{1},\tilde x_{2}, \tilde x_{3})\,d\tilde x_{1}d\tilde x_{2} d\tilde x_{3}\biggr]\sz{e} ^{-i(x_{1}\gamma_{1}+x_{2}\gamma_{2}+x_{3}\gamma_{3})}\,d\gamma_{1}d\gamma_{2}d\gamma_{3} \;= 
\nonumber\\
\nonumber\\
=\; S_{11}(\tilde{f}_{1})\;+\; S_{12}(\tilde{f}_{2})\;+\; S_{13}(\tilde{f}_{3})\;+\;B(u_{1}^0),\;\;\; (\ref{A8})\quad\quad\quad\quad\quad\quad \quad\quad\quad \quad\quad\quad 
\end{eqnarray}

\begin{eqnarray}\label{eqn161}
u_{2}(x_{1}, x_{2}, x_{3}, t)\;=\;
\frac{1}{(2\pi)^{3/2}} \int_{-\infty}^{\infty} \int_{-\infty}^{\infty} \int_{-\infty}^{\infty} \biggl[ \int_{0}^{t} \sz{e} ^{-\nu (\gamma_{1}^{2} +\gamma_{2}^{2} +\gamma_{3}^{2}) (t-\tau)} \frac{[( \gamma_{3}^{2} +\gamma_{1}^{2})  \tilde{F}_{2}( \gamma_{1}, \gamma_{2}, \gamma_{3}, \tau)]} { (\gamma_{1}^{2} +\gamma_{2}^{2}+\gamma_{3}^{2} ) } \,d\tau \;-
\nonumber\\
\nonumber\\
-\; \int_{0}^{t} \sz{e} ^{-\nu (\gamma_{1}^{2} +\gamma_{2}^{2} +\gamma_{3}^{2}) (t-\tau)}\frac{[\gamma_{2}\gamma_{3} \tilde{F}_{3}( \gamma_{1}, \gamma_{2}, \gamma_{3}, \tau ) +\gamma_{2}\gamma_{1} \tilde{F}_{1}( \gamma_{1}, \gamma_{2}, \gamma_{3}, \tau )]} { (\gamma_{1}^{2} +\gamma_{2}^{2}+\gamma_{3}^{2} ) }\,d\tau\;   +
\quad
\nonumber\\
\nonumber\\
+\;\sz{e} ^{-\nu (\gamma_{1}^{2} +\gamma_{2}^{2} +\gamma_{3}^{2}) t} \;U_{2}^{0}(\gamma_{1} ,\gamma_{2} ,\gamma_{3})\biggr] \;\sz{e} ^{-i(x_{1}\gamma_{1}+x_{2}\gamma_{2}+x_{3}\gamma_{3})}\,d\gamma_{1}d\gamma_{2}d\gamma_{3}\;=
\nonumber\\
\nonumber\\
=\;-\;\frac{1}{8\pi^{3}} \int_{-\infty}^{\infty} \int_{-\infty}^{\infty} \int_{-\infty}^{\infty}\frac{ \gamma_{2}\gamma_{1}} { (\gamma_{1}^{2} +\gamma_{2}^{2}+\gamma_{3}^{2} ) }  \biggl[ \int_{0}^{t} \sz{e} ^{-\nu (\gamma_{1}^{2} +\gamma_{2}^{2} +\gamma_{3}^{2}) (t-\tau)} \int_{-\infty}^{\infty}\int_{-\infty}^{\infty}\int_{-\infty}^{\infty}\sz{e}  ^{i(\tilde x_{1}\gamma_{1}+\tilde x_{2}\gamma_{2}+\tilde x_{3}\gamma_{3})} \cdot 
\nonumber\\
\nonumber\\
\cdot \tilde{f}_{1}(\tilde x_{1},\tilde x_{2}, \tilde x_{3},\tau)\,d\tilde x_{1}d\tilde x_{2} d\tilde x_{3}d\tau\biggr]\sz{e} ^{-i(x_{1}\gamma_{1}+x_{2}\gamma_{2}+x_{3}\gamma_{3})}\,d\gamma_{1}d\gamma_{2}d\gamma_{3}\;+
\nonumber\\
\nonumber\\
+\;\frac{1}{8\pi^{3}} \int_{-\infty}^{\infty} \int_{-\infty}^{\infty} \int_{-\infty}^{\infty}\frac{ (\gamma_{3}^2 + \gamma_{1}^2)} { (\gamma_{1}^{2} +\gamma_{2}^{2}+\gamma_{3}^{2} ) }  \biggl[ \int_{0}^{t} \sz{e} ^{-\nu (\gamma_{1}^{2} +\gamma_{2}^{2} +\gamma_{3}^{2}) (t-\tau)} \int_{-\infty}^{\infty}\int_{-\infty}^{\infty}\int_{-\infty}^{\infty}\sz{e}  ^{i(\tilde x_{1}\gamma_{1}+\tilde x_{2}\gamma_{2}+\tilde x_{3}\gamma_{3})} \cdot 
\nonumber\\
\nonumber\\
\cdot \tilde{f}_{2}(\tilde x_{1},\tilde x_{2}, \tilde x_{3},\tau)\,d\tilde x_{1}d\tilde x_{2} d\tilde x_{3}d\tau\biggr]\sz{e} ^{-i(x_{1}\gamma_{1}+x_{2}\gamma_{2}+x_{3}\gamma_{3})}\,d\gamma_{1}d\gamma_{2}d\gamma_{3}\;-
\nonumber\\
\nonumber\\
-\;\frac{1}{8\pi^{3}} \int_{-\infty}^{\infty} \int_{-\infty}^{\infty} \int_{-\infty}^{\infty}\frac{ \gamma_{2}\gamma_{3}} { (\gamma_{1}^{2} +\gamma_{2}^{2}+\gamma_{3}^{2} ) }  \biggl[ \int_{0}^{t} \sz{e} ^{-\nu (\gamma_{1}^{2} +\gamma_{2}^{2} +\gamma_{3}^{2}) (t-\tau)} \int_{-\infty}^{\infty}\int_{-\infty}^{\infty}\int_{-\infty}^{\infty}\sz{e}  ^{i(\tilde x_{1}\gamma_{1}+\tilde x_{2}\gamma_{2}+\tilde x_{3}\gamma_{3})} \cdot 
\nonumber\\
\nonumber\\
\cdot \tilde{f}_{3}(\tilde x_{1},\tilde x_{2}, \tilde x_{3},\tau)\,d\tilde x_{1}d\tilde x_{2} d\tilde x_{3}d\tau\biggr]\sz{e} ^{-i(x_{1}\gamma_{1}+x_{2}\gamma_{2}+x_{3}\gamma_{3})}\,d\gamma_{1}d\gamma_{2}d\gamma_{3}\;+
\nonumber\\
\nonumber\\
+\;\frac{1}{8\pi^{3}} \int_{-\infty}^{\infty} \int_{-\infty}^{\infty} \int_{-\infty}^{\infty} \sz{e} ^{-\nu (\gamma_{1}^{2} +\gamma_{2}^{2} +\gamma_{3}^{2}) t}\biggl[  \int_{-\infty}^{\infty}\int_{-\infty}^{\infty}\int_{-\infty}^{\infty}\sz{e}  ^{i(\tilde x_{1}\gamma_{1}+\tilde x_{2}\gamma_{2}+\tilde x_{3}\gamma_{3})}
 \cdot \quad\quad\quad\quad\quad\quad\quad\quad\quad
\nonumber\\
\nonumber\\
\cdot \; u_{2}^{0}(\tilde x_{1},\tilde x_{2}, \tilde x_{3})\,d\tilde x_{1}d\tilde x_{2} d\tilde x_{3}\biggr]\sz{e} ^{-i(x_{1}\gamma_{1}+x_{2}\gamma_{2}+x_{3}\gamma_{3})}\,d\gamma_{1}d\gamma_{2}d\gamma_{3} \;= 
\nonumber\\
\nonumber\\
=\; S_{21}(\tilde{f}_{1})\;+\; S_{22}(\tilde{f}_{2})\;+\; S_{23}(\tilde{f}_{3})\;+\;B(u_{2}^0),\;\;\; (\ref{A8}) 
\quad\quad\quad\quad\quad\quad \quad\quad\quad \quad\quad\quad 
\end{eqnarray}

\begin{eqnarray}\label{eqn162}
u_{3}(x_{1}, x_{2}, x_{3}, t)\;=\;
\frac{1}{(2\pi)^{3/2}} \int_{-\infty}^{\infty} \int_{-\infty}^{\infty} \int_{-\infty}^{\infty} \biggl[ \int_{0}^{t} \sz{e} ^{-\nu (\gamma_{1}^{2} +\gamma_{2}^{2} +\gamma_{3}^{2}) (t-\tau)} \frac{[( \gamma_{1}^{2} +\gamma_{2}^{2})  \tilde{F}_{3}( \gamma_{1}, \gamma_{2}, \gamma_{3}, \tau)]} { (\gamma_{1}^{2} +\gamma_{2}^{2}+\gamma_{3}^{2} ) } \,d\tau \;-
\nonumber\\
\nonumber\\
-\;\int_{0}^{t} \sz{e} ^{-\nu (\gamma_{1}^{2} +\gamma_{2}^{2} +\gamma_{3}^{2}) (t-\tau)}\frac{[\gamma_{3}\gamma_{1} \tilde{F}_{1}( \gamma_{1}, \gamma_{2}, \gamma_{3}, \tau ) +\gamma_{3}\gamma_{2} \tilde{F}_{2}( \gamma_{1}, \gamma_{2}, \gamma_{3}, \tau )]} { (\gamma_{1}^{2} +\gamma_{2}^{2}+\gamma_{3}^{2} ) }\,d\tau\;   +
\quad
\nonumber\\
\nonumber\\
+\; \sz{e} ^{-\nu (\gamma_{1}^{2} +\gamma_{2}^{2} +\gamma_{3}^{2}) t} \;U_{3}^{0}(\gamma_{1} ,\gamma_{2} ,\gamma_{3})\biggr] \;\sz{e} ^{-i(x_{1}\gamma_{1}+x_{2}\gamma_{2}+x_{3}\gamma_{3})}\,d\gamma_{1}d\gamma_{2}d\gamma_{3}\;=
\nonumber\\
\nonumber\\
=\;-\;\frac{1}{8\pi^{3}} \int_{-\infty}^{\infty} \int_{-\infty}^{\infty} \int_{-\infty}^{\infty}\frac{ \gamma_{3}\gamma_{1}} { (\gamma_{1}^{2} +\gamma_{2}^{2}+\gamma_{3}^{2} ) }  \biggl[ \int_{0}^{t} \sz{e} ^{-\nu (\gamma_{1}^{2} +\gamma_{2}^{2} +\gamma_{3}^{2}) (t-\tau)} \int_{-\infty}^{\infty}\int_{-\infty}^{\infty}\int_{-\infty}^{\infty}\sz{e}  ^{i(\tilde x_{1}\gamma_{1}+\tilde x_{2}\gamma_{2}+\tilde x_{3}\gamma_{3})} \cdot 
\nonumber\\
\nonumber\\
\cdot \tilde{f}_{1}(\tilde x_{1},\tilde x_{2}, \tilde x_{3},\tau)\,d\tilde x_{1}d\tilde x_{2} d\tilde x_{3}d\tau\biggr]\sz{e} ^{-i(x_{1}\gamma_{1}+x_{2}\gamma_{2}+x_{3}\gamma_{3})}\,d\gamma_{1}d\gamma_{2}d\gamma_{3}\;-
\nonumber\\
\nonumber\\
-\;\frac{1}{8\pi^{3}} \int_{-\infty}^{\infty} \int_{-\infty}^{\infty} \int_{-\infty}^{\infty}\frac{ \gamma_{3}\gamma_{2}} { (\gamma_{1}^{2} +\gamma_{2}^{2}+\gamma_{3}^{2} ) }  \biggl[ \int_{0}^{t} \sz{e} ^{-\nu (\gamma_{1}^{2} +\gamma_{2}^{2} +\gamma_{3}^{2}) (t-\tau)} \int_{-\infty}^{\infty}\int_{-\infty}^{\infty}\int_{-\infty}^{\infty}\sz{e}  ^{i(\tilde x_{1}\gamma_{1}+\tilde x_{2}\gamma_{2}+\tilde x_{3}\gamma_{3})} \cdot 
\nonumber\\
\nonumber\\
\cdot \tilde{f}_{2}(\tilde x_{1},\tilde x_{2}, \tilde x_{3},\tau)\,d\tilde x_{1}d\tilde x_{2} d\tilde x_{3}d\tau\biggr]\sz{e} ^{-i(x_{1}\gamma_{1}+x_{2}\gamma_{2}+x_{3}\gamma_{3})}\,d\gamma_{1}d\gamma_{2}d\gamma_{3}\;+
\nonumber\\
\nonumber\\
+\;\frac{1}{8\pi^{3}} \int_{-\infty}^{\infty} \int_{-\infty}^{\infty} \int_{-\infty}^{\infty}\frac{ (\gamma_{1}^2 + \gamma_{2}^2)} { (\gamma_{1}^{2} +\gamma_{2}^{2}+\gamma_{3}^{2} ) }  \biggl[ \int_{0}^{t} \sz{e} ^{-\nu (\gamma_{1}^{2} +\gamma_{2}^{2} +\gamma_{3}^{2}) (t-\tau)} \int_{-\infty}^{\infty}\int_{-\infty}^{\infty}\int_{-\infty}^{\infty}\sz{e}  ^{i(\tilde x_{1}\gamma_{1}+\tilde x_{2}\gamma_{2}+\tilde x_{3}\gamma_{3})} \cdot 
\nonumber\\
\nonumber\\
\cdot \tilde{f}_{3}(\tilde x_{1},\tilde x_{2}, \tilde x_{3},\tau)\,d\tilde x_{1}d\tilde x_{2} d\tilde x_{3}d\tau\biggr]\sz{e} ^{-i(x_{1}\gamma_{1}+x_{2}\gamma_{2}+x_{3}\gamma_{3})}\,d\gamma_{1}d\gamma_{2}d\gamma_{3}\;+
\nonumber\\
\nonumber\\
+\;\frac{1}{8\pi^{3}} \int_{-\infty}^{\infty} \int_{-\infty}^{\infty} \int_{-\infty}^{\infty} \sz{e} ^{-\nu (\gamma_{1}^{2} +\gamma_{2}^{2} +\gamma_{3}^{2}) t}\biggl[  \int_{-\infty}^{\infty}\int_{-\infty}^{\infty}\int_{-\infty}^{\infty}\sz{e}  ^{i(\tilde x_{1}\gamma_{1}+\tilde x_{2}\gamma_{2}+\tilde x_{3}\gamma_{3})}
 \cdot \quad\quad\quad\quad\quad\quad\quad\quad\quad
\nonumber\\
\nonumber\\
\cdot \; u_{3}^{0}(\tilde x_{1},\tilde x_{2}, \tilde x_{3})\,d\tilde x_{1}d\tilde x_{2} d\tilde x_{3}\biggr]\sz{e} ^{-i(x_{1}\gamma_{1}+x_{2}\gamma_{2}+x_{3}\gamma_{3})}\,d\gamma_{1}d\gamma_{2}d\gamma_{3} \;= 
\nonumber\\
\nonumber\\
=\; S_{31}(\tilde{f}_{1})\;+\; S_{32}(\tilde{f}_{2})\;+\; S_{33}(\tilde{f}_{3})\;+\;B(u_{3}^0) ,\;\;\; (\ref{A8})
\quad\quad\quad\quad\quad\quad \quad\quad\quad \quad\quad\quad 
\end{eqnarray}

Here $S_{11}(),\; S_{12}(),\; S_{13}(),\; S_{21}(),\; S_{22}(),\; S_{23}(),\; S_{31}(),\; S_{32}(),\; S_{33}(),\; B()\; $ are the integral operators.

\[S_{12}()\;= \;S_{21}() \]\[S_{13}()\;= \;S_{31}() \]\[S_{23}()\;= \;S_{32}() \]

It follows from $(\ref{eqn160})\;-\; (\ref{eqn162})\;$ that the vector $\vec{u}$ can be represented as:

\begin{equation}\label{eqn164}
\vec{u}\;=\;\bar{\bar{S}}\;\cdot\;\vec{\tilde{f}}\;+\;\bar{\bar{B}}\cdot\vec{u}^{0}\;=\;\bar{\bar{S}}\;\cdot\;\vec{f}\;- \bar{\bar{S}}\cdot(\vec{u}\cdot\nabla)\vec{u}\;+\;\bar{\bar{B}}\cdot\vec{u}^{0}\;,\;\;\;\;\;
\end{equation}

where $\vec{\tilde{f}}$ is determined by formula $(\ref{eqn14})$.

Here $\;\bar{\bar{S}} \; $ and $\bar{\bar{B}}$ are the matrix integral operators:
\[ \left( \begin{array}{ccc}
S_{11} & S_{12} & S_{13} \\
S_{21} & S_{22} & S_{23} \\
S_{31} & S_{32} & S_{33}  \end{array} \right)\]
\[ \left( \begin{array}{ccc}
B & 0 & 0 \\
0 & B & 0 \\
0 & 0 & B  \end{array} \right)\]

\section{Spaces S, $\stackrel{\longrightarrow}{\textbf{TS}}$}$\;\;\cite{GC68}, \cite{RR64}$
$\\$

Let us consider space S of all infinitely differentiable functions $\varphi$(x) defined in N-dimensional space $R^{N}$ (N = 3), such that when $\;\;\mid x \mid \;\rightarrow\; \infty\;\;$ these functions tend to 0, as well as their derivatives of any order, more rapidly than any power of $\frac{1}{\mid x \mid}$.

To define topology in the space S let us introduce countable system of norms 

\begin{equation}\label{eqn200}
\|\varphi\|_{p}\;=\; \mathop{sup_{\;x}}_{\mid k \mid ,\; \mid q \mid \;\leq \;p}\mid x^{k}D^{q} \varphi(x)\mid\;\;\;(p = 0, 1, 2,...)
\end{equation}

where 

\begin{center}
$\mid x^{k}D^{q} \varphi(x)\mid\; = \; \mid x_{1}^{k_{1}}\ldots x_{N}^{k_{N}}\frac{\partial^{q_{1}+\cdots + q_{N}}\varphi(x)}{\partial {x_{1}^{q_{1}}}\ldots \partial {x_{N}^{q_{N}}} }\mid$
\end{center}
\begin{center}
\end{center}
\begin{center}
$k = (k_{1}, \ldots , k_{N}),\;\; q = (q_{1}, \ldots , q_{N}),\;\;x^{k} =  x_{1}^{k_{1}}\ldots x_{N}^{k_{N}} $
\end{center}
\begin{center}
\end{center}
\begin{center}
$D^{q} = \frac{\partial^{q_{1}+\cdots + q_{N}}}{\partial {x_{1}^{q_{1}}}\ldots \partial {x_{N}^{q_{N}}} },\;\; (k_{1}, \ldots , q_{N} = 0, 1, 2, \ldots) $
\end{center}

Space S is a perfect space (complete countably normed space, in which the bounded sets are compact). Space $\overrightarrow{TS}$ of vector-functions $\vec{\varphi}$ is a direct sum of N  perfect spaces S (N = 3) $\cite{VT80}, \cite{RN72}$ :
\nonumber\\

\begin{center}
$\overrightarrow{TS} = S \oplus S \oplus S$.
\end{center}

To define topology in the space $\overrightarrow{TS}$ let us introduce countable system of norms 

\begin{equation}\label{eqn201}
\|\vec{\varphi}\|_{p}\;=\; \sum_{i = 1}^{N}\|\varphi_{i}\|_{p}\; = \sum_{i = 1}^{N}\mathop{sup_{\;x}}_{\mid k \mid ,\; \mid q \mid \;\leq \;p}\mid x^{k}D^{q} \varphi_{i}(x)\mid\;\;\;(p = 0, 1, 2,...),\;\;(N = 3)
\end{equation}

The Fourier transform maps the space S onto the whole space S and maps the space $\overrightarrow{TS}$ onto the whole space $\overrightarrow{TS}$ $\cite{gS01},\;\cite{GC68}$.

\section{Equivalence of Cauchy problem in differential form $\textbf {(\ref{eqn1})}$ - $\textbf {(\ref{eqn3})}$ and in the form of an integral equation}
$\\$

Let us denote solution of the problem $(\ref{eqn1})$ - $(\ref{eqn3})$ as \{$\vec{u}(x_{1}, x_{2}, x_{3}, t)$, p($x_{1}, x_{2}, x_{3},$ t)\}, in other words let us consider infinitely differentiable by t $\in$ [0,$\infty$) vector-function $\vec{u}(x_{1}, x_{2}, x_{3}, t) \in \overrightarrow{TS}$ and infinitely differentiable function p($x_{1}, x_{2}, x_{3},$ t) $\in$ S, that turn equations $(\ref{eqn1})$ , $(\ref{eqn2})$ into identities. Vector-function $\vec{u}(x_{1}, x_{2}, x_{3}, t)$ also satisfies the initial condition $(\ref{eqn3})\;\; (\vec{u}^{0} (x_{1}, x_{2}, x_{3})\in \overrightarrow{TS})$:

\begin{equation}\label{eqn202}
\vec{u}(x_{1}, x_{2}, x_{3}, t)|_{t = 0}\;=\;\vec{u}^{0}(x_{1}, x_{2}, x_{3})
\end{equation}

Let us put \{$\vec{u}(x_{1}, x_{2}, x_{3}, t)$, p($x_{1}, x_{2}, x_{3},$ t)\} into equations $(\ref{eqn1})$ , $(\ref{eqn2})$ and apply Fourier and Laplace transforms to the result identities considering initial condition $(\ref{eqn3})$. After all required operations (as in sections 2 and 3) we receive that vector-function $\vec{u}(x_{1}, x_{2}, x_{3}, t)$ satisfies integral equation:

\begin{equation}\label{eqn203}
\vec{u}\;=\;\bar{\bar{S}}\;\cdot\;\vec{f}\;- \bar{\bar{S}}\cdot(\vec{u}\cdot\nabla)\vec{u}\;+\;\bar{\bar{B}}\cdot\vec{u}^{0}\; = \bar{\bar{S}}^{\nabla}\cdot\vec{u}
\end{equation}

Then vector-function grad p $\in \overrightarrow{TS}$ is defined by equations $(\ref{eqn1})$ where vector-function $\vec{u}$ is received from equation $(\ref{eqn203})$.

Here $\vec{f} \in \overrightarrow{TS}$, $\vec{u}^{0} \in \overrightarrow{TS}$ and $\bar{\bar{S}},\;\bar{\bar{B}},\;\;\bar{\bar{S}}^{\nabla}\;$ are matrix integral operators. 

Vector-functions  
$ \bar{\bar{S}}\cdot\vec{f},\;\;\bar{\bar{B}}\cdot\vec{u}^{0},\;\; \bar{\bar{S}}\cdot(\vec{u}\cdot\nabla)\vec{u}\;$ are also belong $\overrightarrow{TS}$ since Fourier transform maps perfect space $\overrightarrow{TS}$ onto $\overrightarrow{TS}$.

Going from the other side, let us assume that $\vec{u}(x_{1}, x_{2}, x_{3}, t)\in \overrightarrow{TS}$ is continuous in t $\in$ [0,$\infty$) solution of integral equation $(\ref{eqn203})$. Integral-operators $S_{ij}\cdot(\vec{u}\cdot\nabla)\vec{u}$ are continuous in t $\in$ [0,$\infty$) [see $(\ref{eqn160})\;-\;(\ref{eqn162})$]. From here we receive that according to $(\ref{eqn203})$
$\\$
\begin{center}
$\vec{u}(x_{1}, x_{2}, x_{3}, 0)  = \vec{u}^{0}(x_{1}, x_{2}, x_{3})$ 
\nonumber\\
\nonumber\
\end{center}

and also that $\vec{u}(x_{1}, x_{2}, x_{3}, t)$ is differentiable by t $\in$ [0,$\infty$). As described before, the Fourier transform maps perfect space $\overrightarrow{TS}$ on itself. Hence, \{$\vec{u}(x_{1}, x_{2}, x_{3}, t)$ and $p(x_{1}, x_{2}, x_{3}, t)$\} is the solution of the Cauchy problem $(\ref{eqn1})$ - $(\ref{eqn3})$. From here we see that solving the Cauchy problem $(\ref{eqn1})$ - $(\ref{eqn3})$ is equivalent to finding continuous in t $\in$ [0,$\infty$) solution of integral equation $(\ref{eqn203})$.
\nonumber\\

\section{The Caccioppoli-Banach fixed point principle}$\cite{KA64},\;$ $\cite{VT80},\;$ $\cite{WR73},\;$$\; \cite{KS01},\;$$\; \cite{GD03},\;$$\;\cite{ADL97}\;$
$\\$

From the problem statement $\cite{CF06}$ we have $\vec{f}\;\equiv\; 0$. Let us rewrite integral equation $(\ref{eqn203})$ with this condition

\begin{equation}\label{eqn203a}
\vec{u}\;=\;- \bar{\bar{S}}\cdot(\vec{u}\cdot\nabla)\vec{u}\;+\;\bar{\bar{B}}\cdot\vec{u}^{0}
\end{equation}

Let us divide all parts of the integral equation $(\ref{eqn203a})$ by some constant V, that we will define appropriately below. Then we receive modified integral equation equivalent to equation $(\ref{eqn203a})$:

\begin{equation}\label{eqn203ab}
\vec{u}_{\mbox{\tiny V}}\;=\;- \bar{\bar{S}_{\mbox{\tiny V}}}\cdot(\vec{u}_{\mbox{\tiny V}}\cdot\nabla_{\mbox{\tiny V}})\vec{u}_{\mbox{\tiny V}}\;+\;\bar{\bar{B}_{\mbox{\tiny V}}}\cdot\vec{u}_{\mbox{\tiny V}}^{0}\; = \bar{\bar{S}_{\mbox{\tiny V}}}^{\nabla}\cdot\vec{u}_{\mbox{\tiny V}}
\end{equation}

here

\begin{equation}\label{eqn203z}
\vec{u}_{\mbox{\tiny V}}\;=\;\frac{\vec{u}}{V};\;\vec{u}_{\mbox{\tiny V}}^{0}\;=\;\frac{\vec{u}^{0}}{V};\;\nabla_{\mbox{\tiny V}}\;=\;V\cdot\nabla;\;
x_{k\mbox{\tiny V}}\;=\;\frac{x_{k}}{V};\;\gamma_{k\mbox{\tiny V}}\;=\;V\cdot\gamma_{k};\;\nu_{\mbox{\tiny V}}\;=\;\frac{\nu}{V^{2}};\;\;({1\leq k \leq N}) 
\end{equation}


Let us choose the constant V as:
\begin{equation}\label{eqn203e}
max \;\{\|{\vec{u}^0}\|_{C^2}, \;\|{\vec{u}^0}\|_{L_2}, \;(\|{\vec{u}}\|_{C_2} + \|{\vec{u}}\|_{L_2})\}<<\;V 
\end{equation}

Then we have that

\begin{equation}\label{eqn203f}
|{\vec{u}_{\mbox{\tiny V}}}|<\;1\;\;\;,\;\;\;|{\vec{u}^0_{\mbox{\tiny V}}}|<\;1
\end{equation}

Below in property 4 of Matrix integral operator $\bar{\bar{S}_{\mbox{\tiny V}}}^{\nabla}$, and in Theorems 1 - 3, where two vector-functions $\vec{u}_{\mbox{\tiny V}}$  and $\vec{u}_{\mbox{\tiny V}}^{'}$ are presented, constant V is selected as maximum between V($\vec{u}$) and V$^{'}$($\vec{u}^{'}$). 

We can use the fixed point principle to prove existence and uniqueness of the solution of integral equation $(\ref{eqn203ab})$. 

For this purpose we will operate with the following properties of matrix integral operator $\bar{\bar{S}_{\mbox{\tiny V}}}^{\nabla}$:

1.	Matrix integral operator $\bar{\bar{S}_{\mbox{\tiny V}}}^{\nabla}$ depends continuously  on its parameter t $\in$ [0,$\infty$) (based on formulas $(\ref{eqn160})$ - $(\ref{eqn162})$).
$\\$

2.	Matrix integral operator $\bar{\bar{S}_{\mbox{\tiny V}}}^{\nabla}$ maps vector-functions $\vec{u}_{\mbox{\tiny V}}$ from perfect space $\overrightarrow{TS}$ onto perfect space $\overrightarrow{TS}$. This property directly follows from the properties of Fourier transform $\cite{GC68},\;$ and the form of integrands of integral operators $S_{ij},\; B$ (based on formulas $(\ref{eqn160})$ - $(\ref{eqn162})$). 
$\\$

3.	Matrix integral operator $\bar{\bar{S}_{\mbox{\tiny V}}}$ is "quadratic".
$\\$

4. $\|\bar{\bar{S}_{\mbox{\tiny V}}}^{\nabla}\cdot\vec{u}_{\mbox{\tiny V}}\; - \;\bar{\bar{S}_{\mbox{\tiny V}}}^{\nabla}\cdot\vec{u}_{\mbox{\tiny V}}^{'}\|_{p}\;<\;\|\vec{u}_{\mbox{\tiny V}}\;-\;\vec{u}_{\mbox{\tiny V}}^{'} \|_{p}\;$ for any $\vec{u}_{\mbox{\tiny V}},\;\vec{u}_{\mbox{\tiny V}}^{'} \in \overrightarrow{TS}\;\;$ (${\vec{u}_{\mbox{\tiny V}}\;\neq \;\vec{u}_{\mbox{\tiny V}}^{'}})$ and any t $\in$ [0,$\infty$)    
(based on properties 1, 2, 3 and  formulas $(\ref{eqn160})$ - $(\ref{eqn162}),(\ref{eqn203f} )$).
$\\$

Properties mentioned above allow us to prove that matrix integral operator $\bar{\bar{S}_{\mbox{\tiny V}}}^{\nabla}$ is a contraction operator.
\nonumber\\

\textsc{Theorem 1. }\textbf{Contraction operator.}$\;\;\;\textsc{\cite{KA64}}\;$

Matrix integral operator $\bar{\bar{S}_{\mbox{\tiny V}}}^{\nabla}$ maps perfect space $\overrightarrow{TS}$ onto perfect space $\overrightarrow{TS}$, and for any $\vec{u}_{\mbox{\tiny V}},\;\vec{u}_{\mbox{\tiny V}}^{'}\in \overrightarrow{TS}$ 

(${\vec{u}_{\mbox{\tiny V}}\;\neq \;\vec{u}_{\mbox{\tiny V}}^{'}}$)  the condition 4 is valid. 

Then matrix integral operator $\bar{\bar{S}_{\mbox{\tiny V}}}^{\nabla}$ is a contraction operator, i.e. the following condition is true:

\begin{equation}\label{eqn206}
\|\bar{\bar{S}_{\mbox{\tiny V}}}^{\nabla}\cdot\vec{u}_{\mbox{\tiny V}}\; -\;\bar{\bar{S}_{\mbox{\tiny V}}}^{\nabla}\cdot\vec{u}_{\mbox{\tiny V}}^{'}\|_{p}\;\leq\;\alpha\cdot\|\vec{u}_{\mbox{\tiny V}}\;-\;\vec{u}_{\mbox{\tiny V}}^{'} \|_{p}\;
\end{equation}

where $\alpha\;<\;1$ and is independent from $\vec{u}_{\mbox{\tiny V}},\;\vec{u}_{\mbox{\tiny V}}^{'} \in \overrightarrow{TS}$ for any t $\in$ [0,$\infty$).
\nonumber\\

Proof by contradiction.

Let us assume that the opposite is true. Then there exist such $\vec{u}_{\mbox{\tiny V}n},\;\vec{u}^{'}_{\mbox{\tiny V}n} \in \overrightarrow{TS}$ (n=1,2,$\ldots$) and \[\lim_{ n\rightarrow \infty}\vec{u}_{\mbox{\tiny V}n}, \lim_{ n\rightarrow \infty}\vec{u}^{'}_{\mbox{\tiny V}n} \in \overrightarrow{TS}\;\] that 

\begin{equation}\label{eqn207}
\|\bar{\bar{S}_{\mbox{\tiny V}}}^{\nabla}\cdot\vec{u}_{\mbox{\tiny V}n}\; -\;\bar{\bar{S}_{\mbox{\tiny V}}}^{\nabla}\cdot\vec{u}^{'}_{\mbox{\tiny V}n}\|_{p}\;=\;\alpha_{n}\cdot\|\vec{u}_{\mbox{\tiny V}n}\;-\;\vec{u}^{'}_{\mbox{\tiny V}n} \|_{p}\;\;\;\;\;\;\;\;(n=1,2,\ldots;\;\;\alpha_{n}\;\rightarrow\;1 )\;\;
\end{equation}

Then the limiting result in $(\ref{eqn207})$ would lead to equality 
\nonumber\\
\nonumber\
\begin{center}
$\|\bar{\bar{S}_{\mbox{\tiny V}}}^{\nabla}\cdot\vec{u}_{\mbox{\tiny V}}\; - \;\bar{\bar{S}_{\mbox{\tiny V}}}^{\nabla}\cdot\vec{u}_{\mbox{\tiny V}}^{'}\|_{p}\;=\;\|\vec{u}_{\mbox{\tiny V}}\;-\;\vec{u}_{\mbox{\tiny V}}^{'} \|_{p}\;,$
\nonumber\\
\nonumber\
\end{center}

which is contradicting condition 4. Hence, $\bar{\bar{S}_{\mbox{\tiny V}}}^{\nabla}$ is a contraction operator. 
\nonumber\\

\textsc{Theorem 2. }\textbf{Existence and uniqueness of solution.}$\;\;\;\textsc{\cite{KA64}}\;$

Let us consider a contraction operator $\bar{\bar{S}_{\mbox{\tiny V}}}^{\nabla}$. Then there exists the unique solution $\vec{u}_{\mbox{\tiny V}}^{*}$ of equation $(\ref{eqn203ab})$

in
space $\overrightarrow{TS}$ for any t $\in$ [0,$\infty$). Also in this case it is possible to obtain $\vec{u}_{\mbox{\tiny V}}^{*}$ as a limit of sequence  $\{\vec{u}_{\mbox{\tiny V}n}\}$ , 

where

\begin{center}
$\vec{u}_{\mbox{\tiny V},n+1}\; =\;\bar{\bar{S}_{\mbox{\tiny V}}}^{\nabla}(\vec{u}_{\mbox{\tiny V}n}) \;\;\;\;\;\;\;\;(n=0,1,\ldots),$ 
\nonumber\\
\nonumber\
\end{center}

and $\vec{u}_{\mbox{\tiny V}0}\; = \; 0$.

The rate of conversion of the sequence  $\{\vec{u}_{\mbox{\tiny V}n}\}$  to the solution can be defined from the following inequality:

\begin{equation}\label{eqn208}
\|\vec{u}_{\mbox{\tiny V}n}\; -\;\vec{u}_{\mbox{\tiny V}}^{*}\|_{p}\;\leq\;\frac{\alpha^{n}}{(1 - \alpha)}\|\vec{u}_{\mbox{\tiny V}1}\;-\;\vec{u}_{\mbox{\tiny V}0} \|_{p}\;\;\;\;\;\;\;\;\;(n=0,1,\ldots)
\end{equation}

Proof:

It is clear that 

\begin{center}
$\vec{u}_{\mbox{\tiny V},n+1}\; =\;\bar{\bar{S}_{\mbox{\tiny V}}}^{\nabla}(\vec{u}_{\mbox{\tiny V}n}) ,\;\;\;\;\vec{u}_{\mbox{\tiny V}n}\; =\;\bar{\bar{S}_{\mbox{\tiny V}}}^{\nabla}(\vec{u}_{\mbox{\tiny V},n-1})$. 
\nonumber\\
\nonumber\
\end{center}

It follows from $(\ref{eqn206})$ that 
\nonumber\\
\nonumber\
\begin{center}
$\|\vec{u}_{\mbox{\tiny V},n+1}\; -\;\vec{u}_{\mbox{\tiny V}n}\|_{p}\;\leq\;\alpha\cdot\|\vec{u}_{\mbox{\tiny V}n}\;-\;\vec{u}_{\mbox{\tiny V},n-1} \|_{p}$.
\nonumber\\
\nonumber\
\end{center}

Using similar inequalities one after another while decreasing n we will receive: 
\nonumber\\
\nonumber\
\begin{center}
$\|\vec{u}_{\mbox{\tiny V},n+1}\; -\;\vec{u}_{\mbox{\tiny V}n}\|_{p}\;\leq\;\alpha^{n}\cdot\|\vec{u}_{\mbox{\tiny V}1}\;-\;\vec{u}_{\mbox{\tiny V}0} \|_{p}$.
\nonumber\\
\nonumber\
\end{center}

From this result it follows that

\begin{eqnarray}\label{eqn209}
\|\vec{u}_{\mbox{\tiny V},n+l}\; -\;\vec{u}_{\mbox{\tiny V}n}\|_{p}\;\leq\;\|\vec{u}_{\mbox{\tiny V},n+l}\; -\;\vec{u}_{\mbox{\tiny V},n+l-1}\|_{p}+\cdots+\|\vec{u}_{\mbox{\tiny V},n+1}\; -\;\vec{u}_{\mbox{\tiny V}n}\|_{p}\;\leq\;
\nonumber\\
\nonumber\\
\leq(\alpha^{n+l-1}\;+\cdots+\;\alpha^{n})\|\vec{u}_{\mbox{\tiny V}1}\;-\;\vec{u}_{\mbox{\tiny V}0} \|_{p}\leq\frac{\alpha^{n}}{(1 - \alpha)}\|\vec{u}_{\mbox{\tiny V}1}\;-\;\vec{u}_{\mbox{\tiny V}0} \|_{p}.
\quad\quad
\end{eqnarray}

Because of $\alpha^{n} \rightarrow 0$ for n$\;\;\rightarrow\; \infty$ the obtained estimation $(\ref{eqn209})$ shows that sequence  $\{\vec{u}_{\mbox{\tiny V}n}\}$  is a Cauchy 

sequence. Since the space $\overrightarrow{TS}$ is a perfect space, this sequence converges to an element $\vec{u}_{\mbox{\tiny V}}^{*}\in \overrightarrow{TS}$, 

such that $\bar{\bar{S}_{\mbox{\tiny V}}}^{\nabla}(\vec{u}_{\mbox{\tiny V}}^{*})$ has sense. We use inequality $(\ref{eqn206})$ again and have:
\nonumber\\
\nonumber\
\begin{center}
$\|\vec{u}_{\mbox{\tiny V},n+1}\; -\;\bar{\bar{S}_{\mbox{\tiny V}}}^{\nabla}(\vec{u}_{\mbox{\tiny V}}^{*})\|_{p}\;=\;\|\bar{\bar{S}_{\mbox{\tiny V}}}^{\nabla}(\vec{u}_{\mbox{\tiny V}n})\; -\;\bar{\bar{S}_{\mbox{\tiny V}}}^{\nabla}(\vec{u}_{\mbox{\tiny V}}^{*})\|_{p}\;\leq\;\alpha\cdot\|\vec{u}_{\mbox{\tiny V}n}\;-\;\vec{u}_{\mbox{\tiny V}}^{*} \|_{p}\;\;\;\;\;\;\;\;\;(n=0,1,2,\ldots)$
\nonumber\\
\nonumber\
\end{center}

The right part of the above inequality tends to 0 for n$\;\;\rightarrow\; \infty$ and it means that 
$\;\;\;\vec{u}_{\mbox{\tiny V},n+1}\; \rightarrow\;\bar{\bar{S}_{\mbox{\tiny V}}}^{\nabla}(\vec{u}_{\mbox{\tiny V}}^{*})$

and $\;\;\;\vec{u}_{\mbox{\tiny V}}^{*}\; =\;\bar{\bar{S}_{\mbox{\tiny V}}}^{\nabla}(\vec{u}_{\mbox{\tiny V}}^{*})$. In other words, $\;\;\vec{u}_{\mbox{\tiny V}}^{*}\;\;$ is the solution of equation $(\ref{eqn203ab})$.

Uniqueness of the solution also follows from $(\ref{eqn206})$.
In fact, if there would exist another solution $\widetilde{\vec{u}_{\mbox{\tiny V}}} \in \overrightarrow{TS}$, 

then 
\nonumber\\
\nonumber\
\begin{center}
$\|\widetilde{\vec{u}_{\mbox{\tiny V}}}\; -\;\vec{u}_{\mbox{\tiny V}}^{*}\|_{p}\;=\;\|\bar{\bar{S}_{\mbox{\tiny V}}}^{\nabla}(\widetilde{\vec{u}_{\mbox{\tiny V}}})\; -\;\bar{\bar{S}_{\mbox{\tiny V}}}^{\nabla}(\vec{u}_{\mbox{\tiny V}}^{*})\|_{p}\;\leq\;\alpha\cdot\|\widetilde{\vec{u}_{\mbox{\tiny V}}}\;-\;\vec{u}_{\mbox{\tiny V}}^{*} \|_{p}.$
\nonumber\\
\nonumber\
\end{center}

Such situation could happen only if $\|\widetilde{\vec{u}_{\mbox{\tiny V}}}\; -\;\vec{u}_{\mbox{\tiny V}}^{*}\|_{p}\;=\;0$, or $\widetilde{\vec{u}_{\mbox{\tiny V}}}\;=\;\vec{u}_{\mbox{\tiny V}}^{*}$.
\nonumber\\

We can also receive an estimation $(\ref{eqn208})$ from estimation $(\ref{eqn209})$ as a limiting result 
for $l\;\rightarrow\; \infty$.

Now let us show that continuous dependence of operator $\bar{\bar{S}_{\mbox{\tiny V}}}^{\nabla}$ on t leads to continuous dependence of

the solution of the problem on t.

We will say that matrix integral operator $\bar{\bar{S}_{\mbox{\tiny V}}}^{\nabla}$ is continuous in t at a
point $t_{0} \in$ [0,$\infty$), if for any sequence 

$\{t_{n}\} \in$ [0,$\infty$) with $t_{n}\; \rightarrow\; t_{0}$ for $n\;\rightarrow\; \infty$, the following is true:

\begin{equation}\label{eqn210}
\bar{\bar{S}_{\mbox{\tiny V}}}_{t_{n}}^{\nabla}(\vec{u}_{\mbox{\tiny V}})\;\rightarrow\;\bar{\bar{S}_{\mbox{\tiny V}}}_{t_{0}}^{\nabla}(\vec{u}_{\mbox{\tiny V}})\;\;\;\;\;\;\;\;\texttt{for}\;\; \texttt{any} \;\;\;\;\vec{u}_{\mbox{\tiny V}} \in \overrightarrow{TS}.
\end{equation}

From \textsc{Theorem 2} it follows that for any t $\in$ [0,$\infty$) equation $(\ref{eqn203ab})$ has the unique solution, which depends 

on t. Let us denote it as $\vec{u}_{\mbox{\tiny V}t}^{*}$. We will say that solution of equation $(\ref{eqn203ab})$ depends continuously on  $t$ at

$t\;=\;t_{0}$, if for any sequence $\{t_{n}\} \in$ [0,$\infty$) with $t_{n}\; \rightarrow\; t_{0}$ for $n\;\rightarrow\; \infty$, the following is true:
\nonumber\\
\nonumber\
\begin{center}
$\vec{u}_{\mbox{\tiny V}t_{n}}^{*}\;\rightarrow\;\vec{u}_{\mbox{\tiny V}t_{0}}^{*}$.
\nonumber\\
\nonumber\
\end{center}

\textsc{Theorem 3. }\textbf{Continuous dependence of solution on t.}$\;\;\;\textsc{\cite{KA64}}\;$

Let us consider operator $\bar{\bar{S}_{\mbox{\tiny V}}}_{t}^{\nabla}$ that satisfies condition $(\ref{eqn206})$ for any t $\in$ [0,$\infty$), 
where $\alpha$ is independent
 
from t and that operator $\bar{\bar{S}_{\mbox{\tiny V}}}_{t}^{\nabla}$ is continuous in t at a point $t_{0} \in$ [0,$\infty$).
Then for $t\;=\;t_{0}$ the solution of 

equation $(\ref{eqn203ab})$ depends continuously  on t.

Proof: 

Let us consider any t $\in$ [0,$\infty$). We will construct the solution $\vec{u}_{\mbox{\tiny V}t}^{*}$ of equation $(\ref{eqn203ab})$ as a limit of sequence 

 $\{\vec{u}_{\mbox{\tiny V}n}\}$ :

\begin{equation}\label{eqn211}
\vec{u}_{\mbox{\tiny V},n+1}\; =\;\bar{\bar{S}_{\mbox{\tiny V}}}_{t}^{\nabla}(\vec{u}_{\mbox{\tiny V}n}) \;\;\;\;\;\;\;\;(n=0,1,\ldots;\;\;\;\vec{u}_{\mbox{\tiny V}0}\;=\;\vec{u}_{\mbox{\tiny V}t_{0}}^{*})
\end{equation}

Let us rewrite inequality $(\ref{eqn208})$ for n = 0:
\begin{equation}\label{eqn212}
\|\vec{u}_{\mbox{\tiny V}}^{*}\; -\;\vec{u}_{\mbox{\tiny V}0}\|_{p}\;\leq\;\frac{1}{(1 - \alpha)}\|\vec{u}_{\mbox{\tiny V}1}\;-\;\vec{u}_{\mbox{\tiny V}0} \|_{p}
\end{equation}

Since $\vec{u}_{\mbox{\tiny V}t_{0}}^{*}\; =\;\bar{\bar{S}_{\mbox{\tiny V}}}_{t_{0}}^{\nabla}(\vec{u}_{\mbox{\tiny V}t_{0}}^{*})$, then because of $(\ref{eqn211})$ and $(\ref{eqn212})$ we have:

\begin{equation}\label{eqn213}
\|\vec{u}_{\mbox{\tiny V}t}^{*}\; -\;\vec{u}_{\mbox{\tiny V}t_{0}}^{*}\|_{p}\;\leq\;\frac{1}{(1 - \alpha)}\|\vec{u}_{\mbox{\tiny V}1}\;-\;\vec{u}_{\mbox{\tiny V}0} \|_{p}\;=\;\frac{1}{(1 - \alpha)}\|\bar{\bar{S}_{\mbox{\tiny V}}}_{t}^{\nabla}(\vec{u}_{\mbox{\tiny V}t_{0}}^{*})\;-\;\bar{\bar{S}_{\mbox{\tiny V}}}_{t_{0}}^{\nabla}(\vec{u}_{\mbox{\tiny V}t_{0}}^{*}) \|_{p}
\end{equation}

Now with the help of $(\ref{eqn210})$ we obtain the required continuity of $\vec{u}_{\mbox{\tiny V}t}$ for $t\;=\;t_{0}$.
$\\$

Following (\ref{eqn203z}) and (\ref{eqn203e}) we get the result:
\begin{equation}\label{eqn203z1}
\vec{u}\;=\;\vec{u}_{\mbox{\tiny V}}\cdot V,\;\;\;\;\nu\;=\;\nu_{\mbox{\tiny V}}\cdot V^{2}.
\end{equation}

Then vector-function grad p $\in \overrightarrow{TS}$ is defined by equations $(\ref{eqn1})$ where vector-function $\vec{u}$ is received from equation $(\ref{eqn203z1})$. Function p is defined up to an arbitrary constant.

\textbf{In other words there exists the unique set of smooth functions}  $\mathbf{u_{\infty i}(x, t)}$, $\mathbf{p_{\infty}(x, t)}$  \textbf{(i = 1, 2, 3) on} $\mathbf{R^{3} \times [0,\infty)}$ \textbf{that satisfies} $\mathbf{(\ref{eqn1}), (\ref{eqn2}), (\ref{eqn3})}$ \textbf{and}

\begin{equation}\label{eqn186b}
\mathbf{u_{\infty i},\;p_{\infty} \in  C^{\infty}(R^{3} \times [0,\infty)),}
\nonumber\\
\nonumber\\
\end{equation}

Then, using the inequality $\|{\vec{u}}\|_{L_2}\;\leq\;\|{\vec{u}^0}\|_{L_2}\;$from $\;\cite{oL69},\; \cite{LK63}$, we have

\begin{equation}\label{eqn186c}
\mathbf{\int_{R^{3}}|\vec{u}_{\infty}(x, t)|^{2}dx < C } 
\end{equation}

\textbf{for all t} $\mathbf{\geq 0}$.
\nonumber\\

Let us consider $\nu \rightarrow$ 0 in integral operator $\bar{\bar{S}_{\mbox{\tiny V}}}^{\nabla}$.
Then we see that theorems 1-3 are correct also in case of Euler equations, i.e. there exists unique smooth solution in all time range for this case too.

Hence, we can see that when velocity ${\vec{u^0}}\in \overrightarrow{TS}$, the fluid flow is laminar.
Turbulent flow may occur when velocity ${\vec{u^0}}\notin\overrightarrow{TS}$.

\appendix\section*{}
$\\$

The Fourier integral can be stated in the forms: 

\begin{eqnarray}\label{A3}
U( \gamma_{1} , \gamma_{2} , \gamma_{3})=F[\, u(x_{1} , x_{2} , x_{3})]= \frac{1}{(2\pi)^{3/2}} \int_{-\infty}^{\infty} \int_{-\infty}^{\infty} \int_{-\infty}^{\infty} u( x_{1} , x_{2} , x_{3})\,\sz{e}  ^{ i( \gamma_{1} x_{1} + \gamma_{2} x_{2} + \gamma_{3} x_{3}) } dx_{1} dx_{2} dx_{3} 
\nonumber\\
\nonumber\\
u( x_{1} , x_{2} , x_{3})= \frac{1}{(2\pi)^{3/2}} \int_{-\infty}^{\infty} \int_{-\infty}^{\infty} \int_{-\infty}^{\infty} U( \gamma_{1} , \gamma_{2},\gamma_{3}  )\, \sz{e}  ^{- i( \gamma_{1} x_{1} + \gamma_{2} x_{2} + \gamma_{3} x_{3}) } d\gamma_{1} d\gamma_{2} d\gamma_{3} 
\end{eqnarray}
$\\\\$
The Laplace integral is usually stated in the following form:

\begin{equation}\label{A4}
U^{\otimes}(\eta)=L[\,u(t)\,]= \int_{0}^{\infty}u(t)\, \sz{e}  ^{-\eta t}dt
\;\;\;\;\; u(t)=\frac{1}{2\pi i}\int_{c- i \infty }^{c + i \infty} U^{\otimes}(\eta) \,\sz{e}  ^{\eta t}d\eta \;\;\;\;\; c > c_{0}
\end{equation}

\begin{equation}\label{A5}
L[\,u^{'}(t)\,]=\eta \,U^{\otimes}(\eta)-u(0)
\end{equation}

$\\\\$\textsc{The convolution theorem A.1.}$\;\;\;\textsc{{\cite{DP65}}},\;\textsc{{\cite{DW46}}}$
$\\\\$
If integrals
\[ U_{1}^{\otimes}(\eta)= \int_{0}^{\infty}u_{1}(t)\, \sz{e}  ^{-\eta t}d\,t  \;\;\;\;\;\;\;\;\;\; U_{2}^{\otimes}(\eta)= \int_{0}^{\infty}u_{2}(t)\, \sz{e}  ^{-\eta t}d\,t \]

absolutely converge by $Re\, \eta > \sigma_{d}$, then  $U^{\otimes}(\eta)\,= \,U_{1}^{\otimes}(\eta)\, U_{2}^{\otimes}(\eta)$ is Laplace transform of 

\begin{equation}\label{A6}
u(t)=\int_{0}^{t}u_{1}(t-\tau)\,u_{2}(\tau)\,d\,\tau
\\
\end{equation}

Useful \emph{Laplace integral}:

\begin{equation}\label{A7}
L[\,\sz{e}  ^{\eta_{k}t}\,]\,=\,\int_{0}^{\infty}\sz{e}  ^{-(\eta-\eta_{k})\,t}d\,t
\;=\; \frac{1}{(\eta-\eta_{k})}\;\;\;\;\;\;\;\;\;(Re\,\eta\,>\,\eta_{k})
\end{equation}
\nonumber\\

\begin{equation}\label{A8}
\end{equation}
In calculations of integrals $(\ref{eqn160})$ - $(\ref{eqn162})$ for components of velocity $u_{1}, u_{2}, u_{3}$ for the inverse Fourier transforms, we have  each integrand $\tilde{f}(\gamma_{1},\gamma_{2},\gamma_{3}$) as a product of functions $\chi$($\gamma_{1},\gamma_{2},\gamma_{3}$) and $\varphi$($\gamma_{1},\gamma_{2},\gamma_{3}$) , 

\[\tilde{f}(\gamma_{1},\gamma_{2},\gamma_{3}) = \chi(\gamma_{1},\gamma_{2},\gamma_{3})\cdot \varphi(\gamma_{1},\gamma_{2},\gamma_{3}),\;\;\; \]

where $\varphi$($\gamma_{1},\gamma_{2},\gamma_{3}$) belongs to space S (functions of $\gamma_{1},\gamma_{2},\gamma_{3}$)$\cite{GC68}$ and $\chi$($\gamma_{1},\gamma_{2},\gamma_{3}$) is one of the fractions:
$\\$

\begin{center}
$\frac{( \gamma_{2}^{2} +\gamma_{3}^{2})} { (\gamma_{1}^{2} +\gamma_{2}^{2}+\gamma_{3}^{2} ) }$,
$\frac{( \gamma_{1}\cdot\gamma_{2})} { (\gamma_{1}^{2} +\gamma_{2}^{2}+\gamma_{3}^{2} ) }$,
$\frac{( \gamma_{1}\cdot\gamma_{3})} { (\gamma_{1}^{2} +\gamma_{2}^{2}+\gamma_{3}^{2} ) }$,
\end{center}
$\nonumber\\$
\begin{center}
$\frac{( \gamma_{3}^{2} +\gamma_{1}^{2})} { (\gamma_{1}^{2} +\gamma_{2}^{2}+\gamma_{3}^{2} ) }$,
$\frac{( \gamma_{2}\cdot\gamma_{3})} { (\gamma_{1}^{2} +\gamma_{2}^{2}+\gamma_{3}^{2} ) }$,
$\frac{( \gamma_{1}^{2} +\gamma_{2}^{2})} { (\gamma_{1}^{2} +\gamma_{2}^{2}+\gamma_{3}^{2} ) }$
\end{center}
$\nonumber\\$

These fractions are infinitely differentiable functions for $\gamma_{1}\neq 0$, $\gamma_{2}\neq 0$, $\gamma_{3}\neq 0$ with one point of discontinuity $\gamma_{1}$=0, $\gamma_{2}$=0, $\gamma_{3}$=0. (Discontinuities have finite values at this point.)

For these calculations the inverse Fourier transforms are defined as Lebesgue integrals with Cauchy principal values. 
$\\\\$
\textsc{Theorem }: 
$\\\\$
Let us prove that the inverse Fourier transform of $\tilde{f}(\gamma_{1},\gamma_{2},\gamma_{3}) = \chi(\gamma_{1},\gamma_{2},\gamma_{3})\cdot \varphi(\gamma_{1},\gamma_{2},\gamma_{3})$

\[F[\tilde{f}]\equiv \psi(\sigma)\equiv \int_{-\infty}^{\infty} \sz{e}  ^{ i( \gamma,\sigma) }\chi(\gamma)\cdot\varphi(\gamma) d\gamma \;\;\;\;\;\;\;\;\;\;\;\;(A.6.1)\]

$\psi(\sigma)$ as function of $\sigma$ belongs to space S (functions of $\sigma$). I.e., $\psi(\sigma)$ has two criteria:

1) $\psi(\sigma)$ is infinitely differentiable function, 

2) when $\;\;\mid \sigma \mid \;\rightarrow\; \infty\;\;$ $\psi(\sigma)$ tends to 0, as well as its derivatives of any order, more rapidly than any power of $\frac{1}{\mid \sigma \mid}$.
$\\\\$
\textsc{Proof}:

Integral in (A.6.1) admits of differentiation with respect to the parameter $\sigma_{j}$, since the integral obtained after formal differentiation remains absolutely convergent:   

\[\frac{\partial \psi(\sigma)}{\partial \sigma_{j}}\equiv \int_{-\infty}^{\infty} i \cdot \gamma_{j}\cdot\sz{e}  ^{ i( \gamma,\sigma) }\chi(\gamma)\cdot\varphi(\gamma) d\gamma \;\;\;\;\;\;\;\;\;\;\;\;(A.6.2)\]

The properties of function $\varphi(\gamma)$ permit this differentiation to be continued without limit. This means that the function $\psi(\sigma)$ is infinitely differentiable (see criterion 1).

To prove criterion 2 we have created a function with parameter n:

\[f_{n}(\gamma_{1},\gamma_{2},\gamma_{3}) = n\cdot\sz{e}  ^{- \frac{1}{n^2 (\gamma_{1}^{2} +\gamma_{2}^{2}+\gamma_{3}^{2} )}}\cdot\chi(\gamma_{1},\gamma_{2},\gamma_{3})\cdot \varphi(\gamma_{1},\gamma_{2},\gamma_{3})\]

\[\;\;\; n = 1,2,3 ...\]

$f_{n}(\gamma_{1},\gamma_{2},\gamma_{3}$) belongs to space S (functions of $\gamma_{1},\gamma_{2},\gamma_{3}$).
$\\\\$
Then we estimate the integral in formula (A.6.1) using the inverse Fourier transform of $f_{n}(\gamma)$ for $\;\;\mid \sigma \mid \;>> 0,$ $\;\;\mid \sigma \mid \;\rightarrow\; \infty\;\;$:

\[\mid\psi(\sigma)\mid\;\equiv \;\mid\int_{-\infty}^{\infty} \sz{e}  ^{ i( \gamma,\sigma) }\cdot\chi(\gamma)\cdot\varphi(\gamma) d\gamma \mid\;\leq\;\mid n\cdot\int_{-\infty}^{\infty} \sz{e}  ^{ i( \gamma,\sigma) }\cdot\sz{e}  ^{- \frac{1}{n^2 \gamma^{2} }}\cdot\chi(\gamma)\cdot \varphi(\gamma) d\gamma \mid\\\;\;\;\;(A.6.3)\]
\[n >> 1...<\infty\]

Since $f_{n}(\gamma_{1},\gamma_{2},\gamma_{3}$) belongs to space S (functions of $\gamma_{1},\gamma_{2},\gamma_{3}$) then the inverse Fourier transform of
$f_{n}(\gamma_{1},\gamma_{2},\gamma_{3}$) belongs to space S (functions of $\sigma_{1},\sigma_{2},\sigma_{3}$). Hence $F[f_{n}]$ is the function, such that when $\;\;\mid \sigma \mid \;\rightarrow\; \infty\;\;$ this function tends to 0, as well as its derivatives of any order, more rapidly than any power of $\frac{1}{\mid \sigma \mid}$. 

So we got from formula (A.6.3) that $\psi(\sigma)$ is the function, such that when $\;\;\mid \sigma \mid \;\rightarrow\; \infty\;\;$ this function tends to 0, more rapidly than any power of $\frac{1}{\mid \sigma \mid}$. 

Then we estimate the integral in formula (A.6.2). Using formula like (A.6.3) for $\frac{\partial \psi(\sigma)}{\partial \sigma_{j}}$, we obtain that $\frac{\partial \psi(\sigma)}{\partial \sigma_{j}}$  tends to 0 more rapidly than any power of $\frac{1}{\mid \sigma \mid}$. The properties of function $\varphi(\gamma)$ permit this differentiation to be continued without limit and further using of formulas like (A.6.3) to derivatives leads to the conclusion that all  derivatives of function $\psi(\sigma)$ tend to 0 more rapidly than any power of $\frac{1}{\mid \sigma \mid}$. 

We have proved that $\psi(\sigma)$ belongs to space S (functions of $\sigma$).
$\\\\$
\textbf{Acknowledgment}:  We express our sincere gratitude to Professor L. Nirenberg, whose suggestion led to conduction of this research.
We are also very thankful to Professor A.B. Gorstko for helpful friendly discussions.

$\\$

\end{document}